\documentclass[journal,twoside,web]{ieeecolor}
\usepackage{generic}
\usepackage{cite}
\usepackage{amsmath,amssymb,amsfonts}
 \usepackage{algorithm}
\usepackage{algpseudocode}
 \usepackage{subfigure}
\usepackage{graphicx}
 \usepackage{bbm}
\usepackage{hyperref}
\hypersetup{hidelinks=true}
\usepackage{textcomp}
\usepackage{mathrsfs}
\renewcommand{\thesection}{\arabic{section}}
\renewcommand{\thesubsection}{\arabic{section}.\arabic{subsection}}

\newtheorem{theorem}{Theorem}[section]
\newtheorem{remark}{Remark}[section]
\newtheorem{lemma}{Lemma}[section]

\newtheorem{assumption}{Assumption}[section]

\newtheorem{proposition}{Proposition}[section]
\newtheorem{definition}{Definition}[section]

\renewcommand{\descriptionlabel}[1]{\hspace{\labelsep}\textnormal{#1}}

\emergencystretch=\maxdimen
\def\BibTeX{{\rm B\kern-.05em{\sc i\kern-.025em b}\kern-.08em
    T\kern-.1667em\lower.7ex\hbox{E}\kern-.125emX}}
\begin{document}
\title{Primal-dual policy learning for  mean-field stochastic  LQR problem}
\author{Xiushan Jiang,  Dong Wang, Weihai Zhang,  \IEEEmembership{Senior Member, IEEE}, Daniel W. C. Ho, \IEEEmembership{Life Fellow, IEEE}, and Yuanqing Wu,  \IEEEmembership{Senior Member, IEEE}
\thanks{This work was supported in part by  
the National Natural Science Foundation of China under Grants (62573427, 
  62373229),  and in part by the Taishan Scholar Young Expert Program
under Grant (tsqn202408110), and in part by the 
Research Grants Council of Hong Kong Special
Administrative Region, China, under Grants (CityU 11205724, CityU 11206825).
 }
\thanks{Xiushan Jiang is with the College of New Energy, China University of Petroleum (East China), Qingdao, China.  Dong Wang and Weihai Zhang are with the College of Electrical Engineering and
Automation, Shandong University of Science and Technology, Qingdao, China (e-mail: w\_hzhang@163.com). Daniel W. C. Ho is with
the Department of Mathematics, City University of Hong Kong,
	Hong Kong. 
Yuanqing Wu is with the School of Intelligent Systems Engineering,
Sun Yat-sen University, Guangzhou, China.}
}

\maketitle

\begin{abstract}
Integrating data-driven techniques with mechanism-driven insights has recently gained popularity as a powerful learning approach to solving traditional LQR problems for designing intelligent controllers in complex dynamic systems. However, the theoretical understanding of various reinforcement learning  algorithms needs further exploration to enhance their efficiency and safety. In this article, by means of primal-dual optimization tools, we study the  partially model-free design of the mean-field stochastic LQR (MF-SLQR) controller using a policy learning   approach.
Firstly, by designing appropriate optimizing variables, the considered  MF-SLQR  problem is  transformed into a  new static nonconvex constrained  optimization problem  with   equivalence preserved in certain senses. After that, the equivalent formulation of the duality results is constructed via finding the solution of the generalized Lyapunov equation. Then, the strong duality is analyzed, based on which we establish a primal-dual algorithm by   Karush-Kuhn-Tucker   conditions. More importantly, a partially model-free implementation  is also presented, which has a direct connection with the classical policy iteration algorithm. Finally, we use a high-dimensional example to validate our methods.
\end{abstract}

\begin{IEEEkeywords}
 Stochastic linear quadratic problem; mean-field system; reinforcement learning; primal-dual method.
\end{IEEEkeywords}
\section{Introduction}
   In classical control theory, the    Linear Quadratic Regulator (LQR) problem  serves as a foundational framework  for    optimal control.  Initially developed by  Kalman for  deterministic systems  \cite{kalman1960}  and extended by  Wonham to  stochastic system with multiplicative noise \cite{wonham1967}, the LQR framework has been
    systematically advanced through the work of  Anderson and Moore \cite{anderson1989}.
     The optimal   feedback solutions to the LQR or stochastic LQR (SLQR) problems are
      typically obtained through matrix Riccati equations  or semidefinite programming formulations. 
       { 
      Mean field (MF) theory was introduced independently by Lasry and Lions \cite{lions} and by Huang, Caines, and Malham'e \cite{huang2007}, and further developed by Andersson and Djehiche \cite{andersson1} as well as Bensoussan et al. \cite{bensoussan1,bensoussan2}. This framework provides an efficient approach for analyzing large-scale multi-agent systems by modeling interactions among numerous agents, where individual effects are negligible but collective dynamics are significant. In contrast to standard MF games, MF type   control problems focus on optimizing a global performance criterion rather than individual objectives. In particular, MF-type  control within the SLQR framework (referred to as the MF-SLQR problem) has garnered significant attention due to its cooperative structure, which facilitates explicit solution derivation, as illustrated in the works of Yong \cite{yong2013} and Elliott, Li, and Ni \cite{elliott2013}, among others. Moreover, this approach plays a pivotal role across diverse application domains, including economics \cite{guo2024}, micro-grid energy storage networks \cite{book}, and other complex dynamical networks.       } 
 However, traditional MF-SLQR controller design often assume complete knowledge of system dynamics \cite{jiangzhongsu}, which is rarely feasible in real-world applications. Moreover, the coupling effects on MF systems further complicate control design, such that a further study on effective  algorithms   is  necessary.

 Reinforcement learning (RL) has emerged as a prominent data-driven approach for advancing intelligent control in dynamic systems with incomplete environment models, which has successful applications in areas such as robotic control \cite{ruanti}, autonomous driving \cite{pangbo}, and strategic game optimization. In RL, the policy learning (PL)  mechanism iteratively improves control strategies through data collection and policy evaluation,  which coincides with the LQR problem to develop  various advanced control methodologies \cite{gravell2021}.
 However, most  existing PL methods are designed for
 deterministic systems  or
 stochastic systems without MF terms  \cite{oura2024,youbased}.
For MF systems, an online value iteration method for LQR problems with unknown dynamics was proposed in  \cite{bingchang},   focusing  on the continuous-time case. To the best of the authors' knowledge,
relatively fewer results are known for the discrete-time MF stochastic systems.
On the other hand,  a common challenge in PL algorithms for LQR problems lies in their reliance on the persistent excitation (PE) condition {\cite{breschi2022,wang2025,jiangguojing}}, which necessitates the addition of artificial exploration noise during the iterative policy updates. Although this ensures sufficient exploration, it can  affect the system stability and hinder the algorithm convergence,   especially when the  state equation incorporates the mathematical expectation of the process.

Thus, to address these difficulties, this study adopts a primal-dual PL (PD-PL) optimization approach to develop a partially model-free controller design   for
the considered MF stochastic systems. Unlike most existing RL methods, the PD-PL approach utilizes a semidefinite programming (SDP) formulation to establish  a fundamental Lagrangian duality framework for our considered MF-SLQR problem,
 that connects model-free PL algorithms with classical Q-learning techniques \cite{paternain2023,watkins}. In particular,  under specific assumptions about the initial state, the Bellman equation was solved through a data-driven  matrix equation rather than a recursive least-squares algorithm \cite{lee2019}. Hence, this enables faster convergence to the optimal solution, even for a high-dimensional  system. This method has proven effective for applications such as tracking control \cite{ruanti} and stochastic LQR problems with additive noise \cite{liman}. However, dealing with multiplicative noise and MF interactions with the system state brings new challenges.
The main contributions of this paper are as follows:
\begin{itemize}
  \item   In contrast to the fully model-based MF-SLQR problem discussed in \cite{nyh,chenchen,qqy}, the proposed PD-PL framework   relaxes the reliance on accurate system modeling of drift parameters. Building on this, the combination of data-driven and mechanism-based approaches enables reliable control and optimized policies, further advancing the development of the MF-SLQR field.
  \item   We reformulate the MF-SLQR problem as a novel equivalent static optimization problem and prove the strong duality theorem. Thereby, this extends the results of \cite{lee2019,liman}, though the reformulation poses greater challenges in considering the adaptability of the solution.
  \item
Unlike the classic Q-learning algorithm, the novel primal-dual update algorithm is theoretically analyzed for uniform convergence based on a new assumption about the initial state, instead of the traditional PE assumption.
\end{itemize}
In the rest of this paper, we rebuild several equivalent optimization problems from the PD perspective in Section 3, after presenting the considered problem in Section 2. Then main results including the strong duality and partially model-free PD-PL algorithm are given in Section 4. Sections 5 and 6 gives an example and concluding this paper, respectively. 

%

 Notation:
 $\mathbb{C}$: the complex plane; $\mathbb{E}$: the mathematical expectation;
 $\mathcal{R}^l$: the $l$-dimensional Euclidean space;
 $\mathcal{R}^{p\times q}$: the set of all $p\times q$ real matrices;
 $\mathscr{S}_n(\mathscr{S}_n^+, \mathscr{S}_n^{++})$:
 the set of all $n\times n$   real symmetric (symmetric positive semidefinite, positive definite) matrices; $\{\Omega, \mathscr{F},  \mathscr{F}_k, \mathcal{P}\}$: a complete filter probability  space with $\mathscr{F}_k$ being the $\sigma$-algebra
 generated by $\{x_0, w_k, k=0, 1, 2, \cdots, k\}$;
 $\mathbb{N}_+(\mathbb{N})$:
 set  of positive (non-negative) integers;
 $\mathbb{N}_\psi\triangleq $ $\{1, 2, \cdots, 0\}$;
 $l_{\mathscr{F}}^2(\mathbb{N}_0; \mathcal{R}^m)$: a space of $l^2$-summable functions  $\phi_k:    \mathbb{N}_\psi\times \Omega\rightarrow \mathcal{R}^m$;
 $ {\operatorname{diag}}(A, B)$: block diagonal matrix with diagonal blocks $A$ and $B$;
$Tr(M)$: the trace of  a square matrix $M$.

\section{Problem formulation and preliminaries}
\subsection{Preliminaries}
The considered MF-SLQR problem aims to find a close-loop state feedback control to minimize the  {cost function}
 \begin{eqnarray}\label{cost1}
   && {J(z; u)}\nonumber\\
  && \triangleq\sum\limits_{k=0}^\infty
   \mathbb{E}\!\left\{\!\left[\!\begin{array}{cc}
   \!x_k\!\\\!u_k\!\end{array}\right]'\Lambda\left[\begin{array}{cc}
  \! x_k\! \\\! u_k\! \end{array}\right]\! +\!
   \left[\begin{array}{cc}
   \! \mathbb{E}x_k\! \\\! \mathbb{E}u_k\!
   \end{array}\right]'\bar{\Lambda}
   \left[\begin{array}{cc}
    \! \mathbb{E}x_k\! \\ \! \mathbb{E}u_k\!
    \end{array}\right]
   \right\}
\end{eqnarray}
 subject to the following controlled linear discrete-time stochastic system with multiplicative noise
\begin{eqnarray}\label{sys1}
 \left\{\begin{array}{l}
 x_{k+1}=(A_1x_k+\bar{A}_1\mathbb{E}x_k
 +B_1u_k
 +\bar{B}_1\mathbb{E}u_k)\\
 \ \ \ \ \  \ \ \ +(A_2x_k+\bar{A}_2\mathbb{E}x_k+B_2u_k
  +\bar{B}_2\mathbb{E}u_k)w_k,\\
  x_0=z, \ \ k\in\mathbb{N},  
 \end{array}
 \right.
 \end{eqnarray}
 where $\Lambda= {\operatorname{diag}}(Q, R)$ and $\bar{\Lambda}= {\operatorname{diag}}(\bar{Q}, \bar{R})$. 
 represent weight matrices of appropriate dimensions, satisfying the following conditions
 \begin{equation}\label{weightma}
  Q, Q+\bar{Q}\geq 0 ,  \ R, R+\bar{R}> 0.
 \end{equation}
 In (\ref{sys1}),  $\{x_k\}_{k\in {\mathbb{N}} }$, valued in $\mathcal{R}^n$ and $\{u_k\}_{k\in {\mathbb{N}} }$, valued in $\mathcal{R}^m$,
 are the system state  and control input, respectively. The system multiplicative noise $ w_k $ is the independent random variables   with zero mean and variance $1$, defined on the given complete probability space $\{\Omega, \mathscr{F}, \mathscr{F}_k, \mathcal{P}\}$.
System matrices $A_1, A_2, \bar{A}_1, \bar{A}_2\in\mathcal{R}^{n\times n}$ and $B_1, B_2, \bar{B}_1,\bar{B}_2\in\mathcal{R}^{n\times m}$  are given deterministic. {$z\in\mathcal{R}^n$ is the   random  initial state  and is assumed to be independent of the
white noise process. }
For simplicity, the notation  $[A_1, \bar{A}_1; A_2, \bar{A}_2]$ refers to the uncontrolled system (\ref{sys1}) (i.e.,
$u_k\equiv0$). Several new notations are  introduced for the system
matrices:
\begin{eqnarray*}
\ \ \ \ \ &&
\hat{A}_1=A_1+\bar{A}_1,\
\hat{B}_1=B_1+\bar{B}_1,\
\hat{A}_2=A_2+\bar{A}_2,\\
&&\hat{B}_2=B_2+\bar{B}_2, \
\hat{Q}=Q+\bar{Q}, \
\hat{R}=R+\bar{R}.
\end{eqnarray*}
 For the   infinite horizon MF-SLQR analysis \cite{elliott2013,chenchen}, the optimal control  exists and is unique with a linear feedback of the current state. The set of  admissible controls $u_k$ are given as follows:
\begin{eqnarray*}
&&\mathcal{U}_{ad}=\{u\in l_{\mathscr{F}}^2( {\mathbb{N}}; \mathcal{R}^m)
: \{u_k\}_{k\in
 {\mathbb{N}}} \ \ \mbox{is\ a \ mean\ square }\\
&&  \ \ \ \ \ \ \ \ \ \ \ \ \  \mbox{stabilizing\ control\ sequence}\}.
\end{eqnarray*}
 {We adopt the following form of a state-feedback control law for the control policy $u_k$:}  
\begin{equation}\label{conpo}
u_k=Fx_k+\bar{F}\mathbb{E}x_k
\end{equation}
with the  {feedback gain matrices} $F$ and $\bar{F}$,
then the state solution for  system (\ref{sys1}) is denoted as
$x_k^{z; F,\bar{F}}$. 
The standard notions of stability, exact observability and exact detectability are
introduced.
 \begin{definition} 
 {(\cite{chenchen}, Definition 1)}
 System $[A_1, \bar{A}_1; A_2, $ $ \bar{A}_2]$ is called asymptotically   mean-square  stable (ASMS) if for any initial state $ {x_0=z}$, there is $\lim\limits_{k\rightarrow \infty} {\mathbb{E}\|x_k^{z}\|^2}=0$.
System (\ref{sys1}) is called  stabilizable if there exists a  control policy
$u_k\in\mathcal{U}_{ad}$ with feedback form  (\ref{conpo}), such that for any initial state
$ {x_0=z}$, the
closed-loop form of system  (\ref{sys1})
is ASMS, i.e., $
\lim\limits_{k\rightarrow \infty}\mathbb{E}\|x_k^{F, \bar{F}; z}\|^2=0. $
\end{definition}
\begin{definition} 
 {(\cite{chenchen}, Definition 3)}
Consider the system
\begin{eqnarray}\label{sysobser}
\left\{\begin{array}{l}
x_{k+1}=A_1x_k+\bar{A}_1\mathbb{E}x_k+(A_2x_k+\bar{A}_2\mathbb{E}x_k)w_k,\\
y_k={\mathcal{Q}^{1/2}}{\mathbb{X}}_k
\end{array}
\right.
\end{eqnarray}
 with ${\mathcal{Q}}= {\operatorname{diag}}(\hat{Q},Q)$. System (\ref{sysobser}) is said to be exact observable if there exists $T\in\mathbb{N}$ such that
\begin{equation*}
y_k\equiv 0, a.s., \forall k\in\mathbb{N}_T\Rightarrow x_0=0.
\end{equation*}
System (\ref{sysobser}) is said to be exact detectable if
\begin{equation*}
y_k\equiv 0, a.s., \forall k\in\mathbb{N}_T\Rightarrow \lim\limits_{k\rightarrow \infty}\mathbb{E}\|x_k\|^2=0.
\end{equation*}
\end{definition}
{
\begin{remark}
The criteria for system stability, observability, and detectability  are standard results in the theory of optimal control and form the basis for the convergence analysis of the following proposed algorithms. 
To streamline the presentation, the detailed proofs and discussions of these criteria are omitted. Interested readers may refer to our prior works, e.g., \cite{zwh2004,chenchen}, for the complete derivations and further discussions.  
\end{remark}
}
 {
In the following,  we make appropriate modifications to system (\ref{sys1}) by introducing new state variables.}
 
Firstly, set $\mathscr{X}_k=\mathbb{E}[(x_k-\mathbb{E}x_k)
(x_k-\mathbb{E}x_k)']$ and $\bar{\mathscr{X}}_k=
\mathbb{E}x_k(\mathbb{E}x_k)'$, then it follows that
\begin{eqnarray}
\left\{\begin{array}{l}
\mathcal{X}_{k+1}=\mathcal{A}_1\mathcal{X}_{k}
\mathcal{A}_1'+\mathcal{A}_2\mathcal{X}_{k}
\mathcal{A}_2'+\mathcal{A}_3\mathcal{X}_{k}
\mathcal{A}_3',\\
\mathcal{X}_\psi=\left[\begin{array}{cccc}
(\mathbb{E}z)(\mathbb{E}z')& \bf{0}\\ \bf{0} &\mathbb{E}[(z-\mathbb{E}z)
(z-\mathbb{E}z)']
\end{array}
\right],
\end{array}
\right.
\end{eqnarray}
where
\begin{eqnarray*}
&&\mathcal{X}_k=\left[\begin{array}{cccc}
\bar{\mathscr{X}}_k& \bf{0} \\ \bf{0}&\mathscr{X}_k
\end{array}
\right], \
\mathcal{A}_1=\left[\begin{array}{cccc}
A_1+\bar{A}_1&\bf{0} \\ \bf{0}&A_1
\end{array}
\right], \\
&&\mathcal{A}_2=\left[\begin{array}{cccc}
\bf{0} & \bf{0} \\ \bf{0} &A_2
\end{array}
\right], \
\mathcal{A}_3=\left[\begin{array}{cccc}
\bf{0} &\bf{0} \\A_2+\bar{A}_2& \bf{0}
\end{array}
\right].
\end{eqnarray*}
Next, define  $\mathbb{X}_{k}'=\left[\begin{array}{cc}\mathbb{E}x_k&
x_k-\mathbb{E}x_k\end{array}\right]$. Under the state-feedback
control (\ref{conpo}),   system (\ref{sys1}) leads to the system
$[\mathcal{A}_1^{F, \bar{F}}; \mathcal{A}_2^F, \mathcal{A}_3^{F, \bar{F}}]$:
\begin{eqnarray}
 \mathbb{X}_{k+1}=\mathcal{A}_1^{F, \bar{F}}\mathbb{X}_k+\mathcal{A}_2^F
 \mathcal{X}_kw_k+\mathcal{A}_3^{F, \bar{F}}\mathbb{X}_kw_k
  \end{eqnarray}
with initial state $ \mathbb{X}_\psi=\left[\begin{array}{cc}\mathbb{E}z'&
(z-\mathbb{E}z)'\end{array}\right]'$,   where
 \begin{eqnarray}
 &\mathcal{A}_1^{F, \bar{F}}=\left[\begin{array}{ccc}
 \hat{A}_1+\hat{B}_1(F+\bar{F})& \bf{0} \\
  \bf{0} &A_1+B_1F
 \end{array}\right],\label{matha1}\\
 &\mathcal{A}_2^F\!=\! \left[\! \begin{array}{ccc}
 \!\bf{0} \! &\!  \bf{0} \\
 \!\bf{0} \! &  A_2\! +\! B_2F
 \end{array}\! \right],\label{matha2}
  \mathcal{A}_3^{F, \bar{F}}\! =\! \left[\! \begin{array}{ccc}
 \!\bf{0} \! &\!  \bf{0}\! \\
 \! \hat{A}_2\! +\!  \hat{B}_2( \!F\! +\! \bar{F})\! &  \! \bf{0}\!
 \end{array}  \right].\label{matha3}
 \end{eqnarray}
 Based on the operator spectrum   technique \cite{zwh2004}, the set of all stabilizing state-feedback gains of the system (\ref{sys1}) is
defined as
$$
\mathcal{F}:=\{F, \bar{F}\in\mathcal{R}^{m\times n}: \sigma(\mathcal{D}_{\mathcal{A}_1^{F, \bar{F}}, \mathcal{A}_2^{F}, \mathcal{A}_3^{F, \bar{F}}})
\subset \mathcal{D}(0, 1)
\}
$$
with $\mathcal{D}(0, 1):=\{\eta: \eta\in\mathbb{C}, |\eta|<1\}$ and
$\sigma(\mathcal{D}_{\mathcal{A}_1^{F, \bar{F}}, \mathcal{A}_2^{F}, \mathcal{A}_3^{F, \bar{F}}})$
being the spectral set of   the generalized Lyapunov operator
      $\mathcal{D}_{\mathcal{A}_1^{F, \bar{F}}, \mathcal{A}_2^{F}, \mathcal{A}_3^{F, \bar{F}}}$:
      \begin{eqnarray*}
     && \mathcal{D}_{\mathcal{A}_1^{F, \bar{F}}, \mathcal{A}_2^{F}, \mathcal{A}_3^{F, \bar{F}}}:
      \mathcal{X}\in\mathscr{S}_{2n}\\
      &&  \longmapsto  (\mathcal{A}_1^{F, \bar{F}})' \mathcal{X}
      \mathcal{A}_3^{F, \bar{F}}
      +(\mathcal{A}_2^{F})'  \mathcal{X}\mathcal{A}_2^{F} +
 (\mathcal{A}_3^{F, \bar{F}})' \mathcal{X}\mathcal{A}_3^{F, \bar{F}}.
      \end{eqnarray*}

\subsection{  Problem Formulation}
From the standard  MF-SLQR theory \cite{elliott2013,jiang2019}, we know that the minimizer
$ (F^*, \bar{F}^*)=\arg \min_{(F, \bar{F})\in\mathcal{F}}J(z, F, \bar{F})$  is independent of the initial state $z$ which is solved via Riccati equations. For technical reason in designing the PL algorithm, we
 make some modification and  research the following MF-SLQR problem.

{\bf Problem 1} (Modified MF-SLQR Problem).
Given the initial state values $z_l\in\mathcal{R}^n$,
$l\in\{1, 2, \cdots, r\}$, with $\sum_{l=1}^r\mathbb{E}z_l(\mathbb{E}z_l)'=Z_2>0$ and
$\sum_{l=1}^r\mathbb{E}[z_lz_l']=Z_1>Z_2$.
 Find the optimal  state feedback control gains
$(F^*, \bar{F}^*)$, such that
\begin{eqnarray*}
\hat{J}(F^*, \bar{F}^*)&:=&\min\limits_
{(F, \bar{F})\in \mathcal{F}}\hat{J}(z_1, \cdots, z_r; F, \bar{F})\\
&=&\min\limits_
{(F, \bar{F})\in \mathcal{F}}\sum\limits_{l=1}^rJ(z_l; F, \bar{F}),
\end{eqnarray*}
where $J(z_l; F, \bar{F})$ is defined in (\ref{cost1}) with control law (\ref{conpo}) working.
\begin{remark}
Compared to  the traditional MF-SLQR problem, which assumes
a single initial state $x_0\in\mathcal{R}^n$, {\bf Problem 1} introduces a  modification
by considering the cost function with
\begin{equation}
\hat{J}(z_1, \cdots, z_r; F, \bar{F})
=\sum\limits_{l=1}^rJ(z_l; F, \bar{F}).
\end{equation}
According to   standard MF-SLQR theory, the optimal control gains, denoted by
\begin{equation*}
(F^*, \bar{F}^*):={\arg\min}_{(F, \bar{F})\in \mathcal{F}}J(z; F,
\bar{F})
\end{equation*}
can be equivalently found by minimizing the modified cost function:
\begin{equation*}
(F^*, \bar{F}^*)={\arg\min}_{(F, \bar{F})\in \mathcal{F}}\hat{J}
(z_1, \cdots, z_r; F,
\bar{F}).
\end{equation*}
\end{remark}
The following assumptions are used throughout the text.
\begin{assumption}\label{asm1}
Assume that
\begin{itemize}
  \item   The weighting matrices $Q,  \bar{Q} ,  R,  \bar{R} $
  satisfy condition (\ref{weightma}).
  \item   System (\ref{sys1}) is stabilizable;
  \item  System $[A_1, \bar{A}_1;
A_2, \bar{A}_2|\mathcal{Q}]$  is exactly observable or exactly detectable.
\end{itemize}
\end{assumption}
Under Assumption \ref{asm1}, based on the result of Theorem 3 in \cite{qqy},  it is straightforward to conclude that the
   optimal value
  $\hat{J}(F^*, \bar{F}^*)$   for   the MF-SLQR problem  exists  and is
  given by
\begin{eqnarray*}
&&\hat{J}(F^*, \bar{F}^*)=\sum_{l=1}^r [\mathbb{E}(z_l'{P}^*z_l)+
(\mathbb{E}z_l)'({P}^*+\bar{P}^*)(\mathbb{E}z_l)]\\
&& \ \ \ \ \ \ \ \ \ \ \ \ \ \ \ =Tr[Z_1{P}^*+Z_2({P}^*+\bar{P}^*)],
\end{eqnarray*}
where  $\bar{P}=\bar{P}^*$  is  the   unique
 solutions to the following coupled  GAREs,  satisfying $P^*\geq 0$,
 $P^*+\bar{P}^*\geq 0$:
\begin{eqnarray}
P=Q+A_1'PA_1+  A_2'PA_2
-\! [M^{P}]'\! [\Upsilon^{P}]^{-1}\! M^{P},\label{gare1}
\end{eqnarray}
\begin{eqnarray}
\bar{P}=\hat{Q}\!+\! \hat{A}_1'\bar{P}
\hat{A}_1\!+\! \hat{A}_2'   P\hat{A}_2\!
- \! [M^{P, \bar{P}}]'[\Upsilon^{P, \bar{P}}]
^{-1}\! [M^{P, \bar{P}}].\label{gare2}
\end{eqnarray}
Here, the terms $\Upsilon^{P}$, $M^{P}$, $\Upsilon^{P, \bar{P}}$, and
$M^{P, \bar{P}}$ are defined as follows:
\begin{eqnarray}
\Upsilon^{P}&=&R+B_1'PB_1+ B_2'PB_2,
\label{ups1}\\
 M^{P}&=&B_1'PA_1+ B_2'PA_2,\label{m1}\\
 \Upsilon^{P, \bar{P}}&=&\hat{R}+\hat{B}_1'
\bar{P}\hat{B}_1
 + \hat{B}_2'P\hat{B}_2,\label{ups21}\\
 M^{P, \bar{P}}&=&\hat{B}_1'\bar{P}
\hat{A}_1+\hat{B}_2'
P\hat{A}_2.\label{m2}
\end{eqnarray}
The corresponding optimal feedback stabilizing control gains are
\begin{eqnarray}
&&F^*\triangleq-[\Upsilon^{P}]^{-1}M^{P},\label{gain1}\\
&&\bar{F}^*\triangleq-\{[\Upsilon^{P, \bar{P}}]^{-1}M^{P, \bar{P}}
-[\Upsilon^{P}]^{-1}M^{P}\}.\label{gain2}
\end{eqnarray}
 For this MF-SLQR problem, the $Q$-function is defined as
\renewcommand{\arraystretch}{0.8}
\begin{eqnarray}
&&Q^*(x_k, u_k)\nonumber\\
&&:=\mathbb{E}
\left\{\left[\begin{array}{cc}
   x_k-\mathbb{E}x_k\\u_k-\mathbb{E}u_k
   \end{array}\right]'\Lambda\left[\begin{array}{cc}
   x_k-\mathbb{E}x_k\\u_k-\mathbb{E}u_k\end{array}\right]+
   \left[\begin{array}{cc}
   \mathbb{E}x_k\\ \mathbb{E}u_k
   \end{array}\right]\right.\nonumber\\
   &&\ \ \ \ \ \ \
   \left.
  \cdot (\Lambda+\bar{\Lambda})\left[\begin{array}{cc}
    \mathbb{E}x_k\\ \mathbb{E}u_k
    \end{array}\right]
   \right\}+\min_u J(x_{k+1}; u)
\nonumber\\
&&= \mathbb{E}\! \left(\! \left[\! \begin{array}{ccc}x_k-\mathbb{E}x_k
\\u_k-\mathbb{E}u_k
\end{array}\right]'\! \mathfrak{X}^*\! \left[\begin{array}{ccc}
x_k-\mathbb{E}x_k
\\u_k-\mathbb{E}u_k
\end{array}\right]\! \right) \nonumber\\
   &&\ \ \ \ \ \ \  +  \left[\begin{array}{ccc}\mathbb{E}
x_k\\
\mathbb{E}u_k
\end{array}\! \right]'\! \bar{\mathfrak{X}}^*\!
\left[\begin{array}{ccc}\mathbb{E}
x_k\\ \mathbb{E}u_k
\end{array}\right],\nonumber
\end{eqnarray}
\renewcommand{\arraystretch}{1.0}
where
\begin{eqnarray}
&&\mathfrak{X}^*=\left[
\begin{array}{ccc}
  Q+A_1'P^*A_1+A_2'P^*A_2&[M^{(1)}]'
 \\
   M^{(1)}&\Upsilon^{(1)}
\end{array}
\right],\label{qfunction1}\\
&&\bar{\mathfrak{X}}^*=\left[\begin{array}{ccc}
{\small \begin{array}{l}
\hat{Q}+\hat{A}_1'\bar{P}^*
\hat{A}_1+ \hat{A}_2'P^*
\hat{A}_2\end{array}}&[M^{(2)}]'\\
M^{(2)}&\Upsilon^{(2)}
\end{array}\right]\label{qfunction2}
\end{eqnarray}
and  $M^{(i)}, \Upsilon^{(i)}$, $i=1, 2$  are
determined by taking $P=P^*$ and
$\bar{P}=\bar{P}^*$
in equations   (\ref{ups1})-(\ref{m2}). The optimal control input
minimizes the $Q$-function, i.e.,
\begin{equation*}
u_k^*=F^*x_k+\bar{F}^*\mathbb{E}x_k
=\arg\min_{u_k\in\mathcal{U}_{ad}}Q^*(x_k, u_k).
\end{equation*}
Next, we introduce a  PI algorithm to solve the coupled GAREs (\ref{gare1})-(\ref{gare2}), which needs the complete knowledge of the system parameter matrices.
\begin{algorithm}
  \caption{PI Algorithm}
  \begin{algorithmic}[1]  
          \State {\bf Initialization:} Select any stabilizer $(F^0, \bar{F}^0)\in(\mathcal{F}\times \mathcal{F}$). Let   the convergence tolerance
          $\varepsilon >0$  and   $i=0$;
           \State  {\bf Repeat Steps $3-6$;}
            \State  {\bf Policy Evaluation:} Solve $P^i$ and $\bar{P}^i$ from
            the    Lyapunov recursion:
            \begin{eqnarray}\label{pe2}
            &&\left[\begin{array}{cc}
           \bar{P}^i \! & \! \bf{0} \\ \bf{0} \! & \! P^i\end{array}\right]\\
           &&
          \!  = \!
            \tilde{\mathcal{Q}}^{F^i, \bar{F}^i} \!+ \!(\mathcal{A}_1^{F^i, \bar{F}^i})'
             \! \left[\begin{array}{cc}
             \bar{P}^i \! & \! \bf{0} \\ \bf{0} \! & \! P^i\end{array}\right] \!(\mathcal{A}_1^{F^i, \bar{F}^i})\nonumber\\
             && \ \ \ \
            +(\mathcal{A}_3^{F^i, \bar{F}^i})' \left[\begin{array}{cc}
             \bar{P}^i \! & \! \bf{0} \\ \bf{0} \! & \! P^i\end{array}\right](\mathcal{A}_3^{F^i, \bar{F}^i})\nonumber\\
            && \ \ \ \
          \!+ \!(\mathcal{A}_2^{F^i})'  \!\left[\begin{array}{cc}
             \bar{P}^i \! & \! \bf{0} \\ \bf{0} \! & \! P^i\end{array}\right]
              (\mathcal{A}_2^{F^i}).
            \end{eqnarray}
       $\mathcal{A}_1^{F^i, \bar{F}^i}$, $\mathcal{A}_2^{F^i}$, and $\mathcal{A}_3^{F^i, \bar{F}^i}$ are defined as in (\ref{matha1})-(\ref{matha3})  with $F$ and $\bar{F}$ replaced by $F^i$ and $\bar{F}^i$, respectively.  $\tilde{\mathcal{Q}}^{F^i, \bar{F}^i}$ is defined as
        \begin{eqnarray*}
 \tilde{\mathcal{Q}}^{F^i, \bar{F}^i}\!=\! \left[\begin{array}{ccc}
  \! \hat{Q}\! +\! (F^i\! +\! \bar{F}^i)' \hat{R}
 (F^i\! +\! \bar{F}^i)  &  \bf{0} \\
 \bf{0}  \! &\!  Q\! +\! (F^i)'RF^i
 \end{array}\right].  
 \end{eqnarray*}
                         \State  {\bf Policy Improvement:} Update control gains from
               \begin{equation}
             F^{i+1}=-[\Upsilon^{P^i}]^{-1}[M^{P^i}]
             \label{pim1}
             \end{equation}
             and
              \begin{equation}
              \bar{F}^{i+1}=-\{[\Upsilon^{P^i, \bar{P}^i}]^{-1}
             [M^{P^i, \bar{P}^i}]
-[\Upsilon^{P^i}]^{-1}[M^{P^i}]\}.\label{pim2}
             \end{equation}
             $\Upsilon^{P^i}$, $\Upsilon^{P^i, \bar{P}^i}$, $M^{P^i}$,
             and $M^{P^i, \bar{P}^i}$ are
             defined as in (\ref{ups1})-(\ref{m2})  with
              $P$ and $\bar{P}$ replaced by $P^i$ and $\bar{P}^i$, respectively.
             \State  $i\leftarrow i+1$
             \State {\bf Until} $\|F^i-F^{i+1}\|+\|\bar{F}^i-\bar{F}^{i+1}\|\leq \varepsilon$.
    \end{algorithmic}
\end{algorithm}
The convergence of Algorithm 1 has been proved in the following result.
\begin{theorem}\label{th22211}
In the  PI algorithm, the sequences $\{P^i\}_{i=0}^\infty$,
$\{\bar{P}^i\}_{i=0}^\infty$, $\{F^i\}_{i=0}^\infty$, and
$\{\bar{F}^i\}_{i=0}^\infty$ have the following  properties.
\begin{description}
  \item[(1)] For $(F^0, \bar{F}^0)\in(\mathcal{F}\times \mathcal{F})$,  all the control gains $\{(F^i, \bar{F}^i)\}_{i=1}^\infty$   are stabilizing.
  \item[(2)]  $P^*\leq P^{i+1}\leq P^{i}$ and
   $\bar{P}^*\leq \bar{P}^{i+1}\leq \bar{P}^{i}$.
$\lim\limits_{i\rightarrow \infty}\! P^i \!  =\!  P^*$, $\lim\limits_{i\rightarrow \infty}\!  \bar{P}^i\!  =\!  \bar{P}^*$, $\lim\limits_{i\rightarrow \infty}\!  F^i\!  =\!  F^*$, and $\lim\limits_{i\rightarrow \infty}\!  \bar{F}^i\!  =\!  \bar{F}^*$, where
      $(P^*, \bar{P}^*)$ are the solution to GAREs (\ref{gare1})-(\ref{gare2}), and
      ($F^*, \bar{F}^*$) are as given in (\ref{gain1})-(\ref{gain2}).
\end{description}
\end{theorem}
{\textit {Proof}}: The proof is provided in Appendix  A.
\section{Optimization Reformulation from the Primal-Dual Perspective}
In this section, we reformulate the  MF-SLQR Problem  using a  novel  primal-dual optimization perspective.
\subsection{Primal-Dual Problem}
Before  presenting  the alternative formulation through the    primal  and dual problems, we make  some modifications to Problem 1  for  technical reasons.
Introduce the augmented state vector
$$
v_k\triangleq \left[\begin{array}{cccc}x_k'
&u_k'\end{array}\right]'\in\mathcal{R}^{n+m},
$$
then
$$
\mathbb{V}_k\triangleq\left[\begin{array}{cccc}\mathbb{E}x_k'
&\mathbb{E}u_k'&x_k'-\mathbb{E}x_k'
&
u_k'-\mathbb{E}u_k'\end{array}\right]'\in\mathcal{R}^{2n+2m}.
$$
The resulting augmented system  is
\begin{eqnarray}\label{ausys}
 \mathbb{V}_{k+1}&=&\mathscr{A}_1^{F, \bar{F}}\mathbb{V}_k+\mathscr{A}_2^
{F}\mathbb{V}_kw_k+\mathscr{A}_3^
{F, \bar{F}}\mathbb{V}_kw_k,\\
  \mathbb{V}_\psi^z&=&\left[\begin{array}{cccc}
(\mathbb{E}v_\psi)'&(v_\psi-\mathbb{E}v_\psi)'
\end{array}
\right]'\nonumber\\
&=&\mathbb{F}_1 \mathbb{E}z+\mathbb{F}_2( z-\mathbb{E}z)\in\mathcal{R}^{2n+2m}, \nonumber
\end{eqnarray}
where
\[
\mathscr{A}_1^{F, \bar{F}}
=
\begin{bmatrix}
\left[\begin{matrix}
\hat{A}_1 & \hat{B}_1 \\
(F+\bar{F})\hat{A}_1 & (F+\bar{F})\hat{B}_1
\end{matrix}\right] & \mathbf{0} \\
\bf{0} &
\left[\begin{matrix}
A_1 & B_1 \\
F A_1 & F B_1
\end{matrix}\right]
\end{bmatrix},
\]
\[
\mathscr{A}_2^{F}
=
\begin{bmatrix}
\mathbf{0} & \mathbf{0} \\
\mathbf{0} &
\left[\begin{matrix}
A_2 & B_2 \\
F A_2 & F B_2
\end{matrix}\right]
\end{bmatrix},
\
\mathscr{A}_3^{F}
=
\begin{bmatrix}
\mathbf{0} & \mathbf{0} \\
\left[\begin{matrix}
\hat{A}_2 & \hat{B}_2 \\
F \hat{A}_2 & F \hat{B}_2
\end{matrix} \right] & \mathbf{0}
\end{bmatrix}
\]
with
$
\mathbb{F}_1=
\left[\begin{array}{ccccc}
I_n&({F}+\bar{F})'&0&0 \end{array}\right]'
$ and $
\mathbb{F}_2=
\left[\begin{array}{ccccc}
0&0&I_n&{F}' \end{array}\right]'.
$
For brevity, this system is denoted as $[\mathscr{A}_1^{F, \bar{F}};\mathscr{A}_2^{F}, \mathscr{A}_3^{F}].$
The state solution is denoted as
$\mathbb{V}_k^{F, \bar{F}; \mathbb{V}_\psi^z}$.  Obviously, the notation  $\mathbb{V}_k^{F, \bar{F}; \mathbb{V}_\psi^{z_l}}$ refers to  the state solution starting from the  initial state $z_l$.
With respect to the augmented system (\ref{ausys}), the cost function (\ref{cost1})  is redefined as
\begin{align}\label{aucos}
 \tilde{J}(F, \bar{F}; \mathbb{V}_\psi^z)\!\triangleq\!
\mathbb{E} \sum\limits_{k=\psi}^\infty
(\mathbb{V}_k^{F, \bar{F}; \mathbb{V}_\psi^z})'\!
\left[\begin{array}{ccccc}
\Lambda+\bar{\Lambda}& \bf{0}\\
\bf{0} &\Lambda \end{array}\right]\! \mathbb{V}_k^{F, \bar{F}; \mathbb{V}_\psi^z}.
 \end{align}
 Since
 $$
 \mathbb{E}\|\mathbb{V}_k^{F, \bar{F}; \mathbb{V}_\psi^z}\|^2
 =\mathbb{E}\|\mathbb{X}_k^{F, \bar{F}; z}\|^2+\mathbb{E}\left\|\left[\begin{array}{ccccc}
F+\bar{F}& \bf{0} \\
\bf{0} &F \end{array}\right]\mathbb{X}_k^{F, \bar{F}; z}\right\|^2,
 $$
 we have
 $$
 \lim\limits_{k\rightarrow \infty}\mathbb{E}\|\mathbb{X}_k^{F, \bar{F}; z}\|^2=0 \Leftrightarrow\lim\limits_{k\rightarrow \infty}\mathbb{E}
 \|\mathbb{V}_k^{F, \bar{F}; \mathbb{V}_0^z}\|^2=0,
 $$
 and the following result is  directly obtained  from   {Theorem 1 of \cite{chenchen}}.
 \begin{lemma}\label{lemm213}
 $(F, \bar{F})\in \mathcal{F}$ iff the augmented system $[\mathscr{A}_1^{F, \bar{F}};\mathscr{A}_2^{F}, \mathscr{A}_3^{F}]$ is
 ASMS.
 \end{lemma}
 The static   optimization problems are given as follows, which can be proved to equivalent with the considered  Problem 1 in some certain senses.

{\bf  Problem 2} (Primal  Problem  I).
Solve the following non-convex minimization problem with variables $\tilde{\mathbbm{S}}_{\mathcal{P}_{I}}\in\mathscr{S}_{2n+2m}$ and  $(F_{\mathcal{P}_{I}}\times\bar{F}_{\mathcal{P}_{I}})\in (\mathcal{R}^{m \times n}\times
\mathcal{R}^{m \times n})$:
\begin{eqnarray}\label{gle1}
\left\{\begin{array}{l}
J_{\mathcal{P}_I}\\
\! \triangleq \!
\inf\limits_{\tilde{\mathbbm{S}}_{\mathcal{P}_{I}}\in\mathscr{S}_{2n+2m}, F_{\mathcal{P}_{I}}, \! \bar{F}_{\mathcal{P}_{I}}\in\mathcal{R}^{m\times n}}
\! Tr\! \Big(\! \left[\! \begin{array}{ccccc}
\! \Lambda+\bar{\Lambda} \! & \! \bf{0} \\
  \! \bf{0} \! &\! \Lambda \end{array}\right]\tilde{\mathbbm{S}}_{\mathcal{P}_{I}}\! \Big),\\
\text{s.t.}\\
 \ \ \ \ \mathscr{A}_1^{F, \bar{F}}\tilde{\mathbbm{S}}_{\mathcal{P}_{I}} (\mathscr{A}_1^{F, \bar{F}})'
+\mathscr{A}_2^{F}\tilde{\mathbbm{S}}_{\mathcal{P}_{I}} (\mathscr{A}_2^{F})' \!
+\! \mathscr{A}_3^{F}\\
\ \ \ \   \cdot  \tilde{\mathbbm{S}}_{\mathcal{P}_{I}}(\mathscr{A}_3^{F})'+\mathbb{F}_1
Z_2\mathbb{F}_1'
 +\mathbb{F}_2(Z_1-Z_2)\mathbb{F}_2'=\tilde{\mathbbm{S}}_{\mathcal{P}_{I}},\\
\ \ \ \ \
(F_{\mathcal{P}_{I}}, \bar{F}_{\mathcal{P}_{I}})\in\mathcal{F}.
 \end{array}
\right.
\end{eqnarray}
{\bf  Problem 3} (Primal  Problem  II).
Solve the following non-convex minimization problem $\tilde{\mathfrak{X}}_{\mathcal{P}_{II}}\in\mathscr{S}_{2n+2m}$ and
$(F_{\mathcal{P}_{II}}\times\bar{F}_{\mathcal{P}_{II}})\in (\mathcal{R}^{m \times n}\times
\mathcal{R}^{m \times n})$:
\begin{eqnarray}\label{gle2}
\left\{\begin{array}{l}
J_{\mathcal{P}_{II}}\triangleq
\inf\limits_{\tilde{\mathfrak{X}}_{\mathcal{P}_{II}}\in\mathscr{S}_{2n+2m}, F_{\mathcal{P}_{II}}, \bar{F}_{\mathcal{P}_{II}}\in\mathcal{R}^{m\times n}}
Tr  (\hbar) ,\\
\text{s.t.}\
 (\mathscr{A}_1^{F, \bar{F}})'\tilde{\mathfrak{X}}_{\mathcal{P}_{II}}\mathscr{A}_1^{F, \bar{F}}
+(\mathscr{A}_2^{F})'\tilde{\mathfrak{X}}
_{\mathcal{P}_{II}}\mathscr{A}_2^{F}
\\
\ \ \ \ \ \   +(\mathscr{A}_3^{F})'\tilde{\mathfrak{X}}_{\mathcal{P}_{II}}
\mathscr{A}_3^{F}+\left[\begin{array}{ccccc}
\Lambda+\bar{\Lambda}& \bf{0} \\
\bf{0} &\Lambda \end{array}\right]=\tilde{\mathfrak{X}}_{\mathcal{P}_{II}},\\
\ \ \ \ \
(F_{\mathcal{P}_{II}}, \bar{F}_{\mathcal{P}_{II}})\in\mathcal{F}.
 \end{array}
\right.
\end{eqnarray}
with
$\hbar= \left(\mathbb{F}_1
Z_2\mathbb{F}_1'+
\mathbb{F}_2(Z_1-Z_2)
\mathbb{F}_2'\right)\tilde{\mathfrak{X}}_{\mathcal{P}_{II}}.
$

{\bf  Problem 4} (Dual  Problem).
Solve
\begin{eqnarray*}
J_{\mathcal{D}}&\triangleq&
\sup\limits_{\tilde{\mathfrak{X}}\in\mathscr{S}_{2n+2m}}d(\tilde{\mathfrak{X}})\\
&=& \sup\limits_{\tilde{\mathfrak{X}}\in\mathscr{S}_{2n+2m}}
\inf\limits_{\tilde{\mathbbm{S}}\in\mathscr{S}_{2n+2m}^+, F, \bar{F}\in\mathcal{F}}
L(\tilde{\mathfrak{X}}, F, \bar{F}, \tilde{\mathbbm{S}}),
\end{eqnarray*}
where
\begin{eqnarray*}
&&L(\tilde{\mathfrak{X}}, F, \bar{F}, \tilde{\mathbbm{S}})\\
&&=Tr\Big(\left[\begin{array}{ccccc}
\Lambda+\bar{\Lambda}& \bf{0} \\
\bf{0} &\Lambda \end{array}\right]\tilde{\mathbbm{S}}\Big)
+Tr\Big\{\Big(\mathscr{A}_1^{F, \bar{F}}\tilde{\mathbbm{S}} (\mathscr{A}_1^{F, \bar{F}})' \\
&&\ \ \ \
+\mathscr{A}_2^{F}\tilde{\mathbbm{S}}(\mathscr{A}_2^{F})'
+\mathscr{A}_3^{F}\tilde{\mathbbm{S}} (\mathscr{A}_3^{F})'
-\tilde{\mathbbm{S}}\\
&&\ \ \ \
 +\mathbb{F}_1
Z_2\mathbb{F}_1'
 +\mathbb{F}_2(Z_1-Z_2)\mathbb{F}_2'\Big)\tilde{\mathfrak{X}}\Big\}.
\end{eqnarray*}
Note that $d(\tilde{\mathfrak{X}})\triangleq \inf\limits_{\tilde{\mathbbm{S}}\in\mathscr{S}_{2n+2m}^+, F, \bar{F}\in\mathcal{F}}
L(\tilde{\mathfrak{X}}, F, \bar{F}, \tilde{\mathbbm{S}})$ is the Lagrangian dual function.
\begin{remark}
 Define the  adjoint operator of $\mathcal{D}_{\mathscr{A}_1^{F, \bar{F}},\mathscr{A}_2^{F}, \mathscr{A}_3^{F}}$  as $\mathcal{D}^*_{\mathscr{A}_1^{F, \bar{F}},\mathscr{A}_2^{F}, \mathscr{A}_3^{F}}$, where
\begin{eqnarray*}
&&\mathcal{D}^*_{\mathscr{A}_1^{F, \bar{F}},\mathscr{A}_2^{F}, \mathscr{A}_3^{F}}:  \tilde{\mathbbm{S}} \in\mathscr{S}_{2n+2m} \\
&& \longmapsto \mathscr{A}_1^{F, \bar{F}}\tilde{\mathbbm{S}}  (\mathscr{A}_1^{F, \bar{F}})'
+\mathscr{A}_2^{F}\tilde{\mathbbm{S}}  (\mathscr{A}_2^{F})' \!
+\! \mathscr{A}_3^{F} \tilde{\mathbbm{S}} (\mathscr{A}_3^{F})'.
\end{eqnarray*}
 Since all  system parameter matrices are real,  it follows that $$\sigma(\mathcal{D}^*_{\mathscr{A}_1^{F, \bar{F}},\mathscr{A}_2^{F}, \mathscr{A}_3^{F}})=\sigma(\mathcal{D}_{\mathscr{A}_1^{F, \bar{F}},\mathscr{A}_2^{F}, \mathscr{A}_3^{F}}).
 $$
According to Lemma \ref{lemm213},  the condition
$(F, \bar{F})\in\mathcal{F}$  is  equivalent to   $[\mathscr{A}_1^{F, \bar{F}};\mathscr{A}_2^{F}, \mathscr{A}_3^{F}]$  being  ASMS, which  implies that
 $\sigma(\mathcal{D}_{\mathscr{A}_1^{F, \bar{F}},\mathscr{A}_2^{F}, \mathscr{A}_3^{F}})\subset \mathcal{D}(0, 1)$.  Therefore,  $(F, \bar{F})\in\mathcal{F}$ implies that   GLEs (\ref{gle1})-(\ref{gle2}) have positive  semidefinite solutions.
\end{remark}

In Primal Problems  I and  II, the constraint $(F, \bar{F}) \in \mathcal{F}$ represents a non-convex condition. The following result shows that, by utilizing the exact detectability and the stochastic PBH eigenvector criteria, this condition can be replaced with $\tilde{\mathbbm{S}} \geq 0$, thereby we are able to transform  Problem  I into a convex formulation.

\begin{proposition}\label{proposition123}
 In Primal Problem  I, the constrained condition $(F, \bar{F})\in \mathcal{F}$ can be replaced by $\tilde{\mathbbm{S}}_{\mathcal{P}_{I}}\geq 0$.
\end{proposition}
{\textit{Proof:}} The proof is provided in Appendix B.
\subsection{Strong Duality}
The dynamic  MF-SLQR optimization problem (Problem 1) has been  reformulated into the constrained static primal-dual optimization problems (Problems 2-4), with their relationships illustrated in Fig. 1.
Before presenting the strong duality analysis, we first demonstrate that Primal Problems  I and  II are the equivalent transformations.
 \begin{proposition}\label{proposition1}
Both Primal Problem  I and Primal Problem  II have unique optimal solutions, denoted by \((\tilde{\mathbbm{S}}_{\mathcal{P}_I}^*, F_{\mathcal{P}_I}^*, \bar{F}_{\mathcal{P}_I}^*)\) and \((\tilde{\mathfrak{X}}_{\mathcal{P}_{II}}^*, F_{\mathcal{P}_{II}}^*, \bar{F}_{\mathcal{P}_{II}}^*)\), respectively. Furthermore, both Primal Problem  I and Primal Problem  II are equivalent to Problem 1, in the sense that \(J_{\mathcal{P}_I} =J_{\mathcal{P}_{II}}= J(F^*, \bar{F}^*)\), \(F_{\mathcal{P}_I}^* = F^*\), and \(\bar{F}_{\mathcal{P}_I}^*=\bar{F}_{\mathcal{P}_{II}}^*
= \bar{F}^*\).
\end{proposition}
{\textit{Proof:}} The proof is provided in Appendix C.\\
The solution  matrix $\tilde{\mathbbm{S}}_{\mathcal{P}_I}\triangleq
\left[\begin{array}{ccccc}
\bar{\mathbbm{S}}_{\mathcal{P}_I} & \bf{0} \\
 \bf{0} & \mathbbm{S}_{\mathcal{P}_I}
 \end{array}\right]$  of GLE (\ref{gle1})  has the following
 forms:
 $$\bar{\mathbbm{S}}_{\mathcal{P}_I}\triangleq
\left[\begin{array}{ccccc}
{\bar{\mathbbm{S}}}_{\mathcal{P}_I}^{11}&{\bar{\mathbbm{S}}}
_{\mathcal{P}_I}^{12}\\
{\bar{\mathbbm{S}}}_{\mathcal{P}_I}^{12'}&{
\bar{\mathbbm{S}}}_{\mathcal{P}_I}^{22}
 \end{array}\right], \  {\mathbbm{S}}_{\mathcal{P}_I}\triangleq
\left[\begin{array}{ccccc}
 {\mathbbm{S}}_{\mathcal{P}_I}^{11}& {\mathbbm{S}}_{\mathcal{P}_I}^{12}\\
 {\mathbbm{S}}_{\mathcal{P}_I}^{12'}&{\mathbbm{S}}_{\mathcal{P}_I}^{22} \end{array}\right]
  $$
  with  $({\bar{\mathbbm{S}}}_{\mathcal{P}_I}^{11}\times {\bar{\mathbbm{S}}}_{\mathcal{P}_I}^{12}\times {\bar{\mathbbm{S}}}_{\mathcal{P}_I}^{22})\in(\mathscr{S}_{n}\times
  \mathcal{R}^{n\times m}\times \mathscr{S}_m)$ and $({ {\mathbbm{S}}}_{\mathcal{P}_I}^{11}\times { {\mathbbm{S}}}_{\mathcal{P}_I}^{12}\times { {\mathbbm{S}}}_{\mathcal{P}_I}^{22})\in(\mathscr{S}_{n}\times
  \mathcal{R}^{n\times m}\times \mathscr{S}_m)$.
  The
   properties for  the solution $\tilde{\mathbbm{S}}_{\mathcal{P}_I}$ are given in the following lemma.
\begin{lemma}\label{lemmaproblem2}
In Primal Problem  I, any feasible solution $\tilde{\mathbbm{S}}\in\mathscr{S}_{2n+2m}$
and $(F\times\bar{F})\in (\mathcal{R}^{m \times n}\times
\mathcal{R}^{m \times n})$ satisfy
\begin{enumerate}
  \item  \begin{align*}
  &\left[\begin{array}{ccccc}
{\bar{\mathbbm{S}}}_{\mathcal{P}_I}^{11}& \bf{0} \\ \bf{0}
& { {\mathbbm{S}}}_{\mathcal{P}_I}^{11}\end{array}\right] \\
&  =
\mathcal{A}_1^{F, \bar{F}} \left[\begin{array}{ccccc}
{\bar{\mathbbm{S}}}_{\mathcal{P}_I}^{11}& \bf{0} \\ \bf{0}
& { {\mathbbm{S}}}_{\mathcal{P}_I}^{11}\end{array}\right]
(\mathcal{A}_1^{F, \bar{F}})'+
 \left[\begin{array}{ccccc}
Z_2& \bf{0} \\ \bf{0}
& Z_1-Z_2\end{array}\right]\\
&  \!+\!
\mathcal{A}_2^{F} \!\left[\begin{array}{ccccc}\!
{\bar{\mathbbm{S}}}_{\mathcal{P}_I}^{11}\! & \! \bf{0}\!  \\ \! \bf{0}\!
& \! { {\mathbbm{S}}}_{\mathcal{P}_I}^{11}\! \end{array}\right]\!
(\mathcal{A}_2^{F})'
\! +\!
\mathcal{A}_3^{F, \bar{F}}\! \left[\begin{array}{ccccc}
\! {\bar{\mathbbm{S}}}_{\mathcal{P}_I}^{11}\! & \! \bf{0} \! \\ \! \bf{0}\!
& \! { {\mathbbm{S}}}_{\mathcal{P}_I}^{11}\!\end{array}\right]\!
(\mathcal{A}_3^{F, \bar{F}})'\! .
\end{align*}
  \item  \begin{align}
&F={ {\mathbbm{S}}}_{\mathcal{P}_I}^{12'}({ {\mathbbm{S}}}_{\mathcal{P}_I}^{11})^{-1},\label{fsfs1}\\
&\bar{F}={\bar{\mathbbm{S}}}_{\mathcal{P}_I}^{12'}
({\bar{\mathbbm{S}}}_{\mathcal{P}_I}^{11})^{-1}
-{ {\mathbbm{S}}}_{\mathcal{P}_I}^{12'}({ {\mathbbm{S}}}_{\mathcal{P}_I}^{11})^{-1}\label{fsfs2}.
\end{align}
\end{enumerate}
 \end{lemma}
The solution matrix $\tilde{\mathfrak{X}}$ of GLE (\ref{gle2}) is uniquely expressed as $\tilde{\mathfrak{X}}=\sum_{k=0}^\infty \mathscr{Y}_k$, where  $\mathscr{Y}_k$ solves
\begin{eqnarray*}
 \left\{\begin{array}{l}
 \mathscr{Y}_{k+1}= (\mathscr{A}_1^{F, \bar{F}})'\mathscr{Y}_k\mathscr{A}_1^{F, \bar{F}}
+(\mathscr{A}_2^{F})'\mathscr{Y}_k\mathscr{A}_2^{F}\\
\ \ \ \ \ \ \ \  \ \ \
+(\mathscr{A}_3^{F})'\mathscr{Y}_k  \mathscr{A}_3^{F},\\
 \mathscr{Y}_0=\left[\begin{array}{ccccc}
\Lambda+\bar{\Lambda}& \bf{0} \\
 \bf{0} &\Lambda \end{array}\right].
 \end{array}\right.
\end{eqnarray*}
Similar to the properties of feasible solutions presented in Lemma \ref{lemmaproblem2}, the properties of the optimal solution for Primal Problem II are outlined below.
\begin{lemma}
For Primal Problem  II, the optimal solution is $(\tilde{\mathfrak{X}}_{\mathcal{P}_{II}}^*, F_{\mathcal{P}_{II}^*}, \bar{F}_{\mathcal{P}_{II}^*})$,  where
\begin{eqnarray*}
&&\tilde{\mathfrak{X}}_{\mathcal{P}_{II}}^*
       \triangleq   \left[\begin{array}{ccccc}
\bar{\mathfrak{X}}_{\mathcal{P}_{II}}^*& \bf{0} \\
 \bf{0} &\mathfrak{X}_{\mathcal{P}_{II}}^*\end{array}\right], \label{xllaaa1}\\
&&\bar{\mathfrak{X}}_{\mathcal{P}_{II}}^*
       \triangleq   \left[\begin{array}{ccccc}
\bar{\mathfrak{X}}_{\mathcal{P}_{II}}^{11*}&
\bar{\mathfrak{X}}_{\mathcal{P}_{II}}^{12*}\\
*& \bar{\mathfrak{X}}_{\mathcal{P}_{II}}^{22*}\end{array}\right], \\
  &&{\mathfrak{X}}_{\mathcal{P}_{II}}^*
       \triangleq   \left[\begin{array}{ccccc}
 {\mathfrak{X}}_{\mathcal{P}_{II}}^{11*}&
 {\mathfrak{X}}_{\mathcal{P}_{II}}^{12*}\\
*& {\mathfrak{X}}_{\mathcal{P}_{II}}^{22*}\end{array}\right]\label{xllaaa2}
\end{eqnarray*}
with  $(\bar{\mathfrak{X}}_{\mathcal{P}_{II}}^*, \mathfrak{X}_{\mathcal{P}_{II}}^*)\in(\mathscr{S}_{n+m}\times
\mathscr{S}_{n+m})$,
$(\bar{\mathfrak{X}}_{\mathcal{P}_{II}}^{11*}\times
\bar{\mathfrak{X}}_{\mathcal{P}_{II}}^{12*}\times
\bar{\mathfrak{X}}_{\mathcal{P}_{II}}^{22*})
\in(\mathscr{S}_{n}\times
  \mathcal{R}^{n\times m}\times \mathscr{S}_m)$,  and $({\mathfrak{X}}_{\mathcal{P}_{II}}^{11*}\times
{\mathfrak{X}}_{\mathcal{P}_{II}}^{12*}\times
{\mathfrak{X}}_{\mathcal{P}_{II}}^{22*})
\in(\mathscr{S}_{n}\times
  \mathcal{R}^{n\times m}\times \mathscr{S}_m)$,   satisfies
 \begin{enumerate}
   \item   $\bar{\mathfrak{X}}_{\mathcal{P}_{II}}^*
=\bar{\mathfrak{X}}^*$ and $\mathfrak{X}_{\mathcal{P}_{II}}^*
=\mathfrak{X}^*$ with
$\bar{\mathfrak{X}}^*$
 and ${\mathfrak{X}}^*$ being the $Q$-function matrices defined in
 (\ref{qfunction1})-(\ref{qfunction2}).
   \item $\tilde{\mathfrak{X}}_{\mathcal{P}_{II}}^*\geq 0$,
     $\bar{\mathfrak{X}}_{\mathcal{P}_{II}}^{22*}>0$, ${\mathfrak{X}}_{\mathcal{P}_{II}}^{22*}>0$,
     \begin{eqnarray*}
     &&F_{\mathcal{P}_{II}}^*=
     ({\mathfrak{X}}_{\mathcal{P}_{II}}^{22*})^{-1}
     ({\mathfrak{X}}_{\mathcal{P}_{II}}^{12*})', \\
     &&\bar{F}_{\mathcal{P}_{II}}^*= (\bar{\mathfrak{X}}_{\mathcal{P}_{II}}^{22*})^{-1}
     (\bar{\mathfrak{X}}_{\mathcal{P}_{II}}^{12*})'-
      ({\mathfrak{X}}_{\mathcal{P}_{II}}^{22*})^{-1}
     ({\mathfrak{X}}_{\mathcal{P}_{II}}^{12*})'.
     \end{eqnarray*}
 \end{enumerate}
\end{lemma}
{\textit{Proof:}} The proof is provided in Appendix D. \\
The weak duality theorem in \cite{weakduality} demonstrates that $J_{\mathcal{D}} \leq J_{\mathcal{P}_I}$ holds naturally. Next, we prove that strong duality also holds.
\begin{theorem}\label{sdth}
(Strong Duality) $J_{\mathcal{D}} = J_{\mathcal{P}_I}$
\end{theorem}
{\textit{Proof:}} The proof is provided in Appendix E.
\section{Primal-Dual Adaptive MF-SLQR Design}
In this section, we shall find a policy learning algorithm with that partial information in system (\ref{sys1}) is unknown.
For this, we first made a modification for Primal Problems I-II. \\
For the augmented state vector $v_k\triangleq \left[\begin{array}{ccc}
x_k'&u_k'\end{array}\right]'$, assume the initial state $v_\psi=\left[\begin{array}{ccc}
x_\psi'&u_\psi'\end{array}\right]'$  is known in advance, where the initial controller $u_\psi$ is freely chosen. This implies that the control policy $u_k=Fx_k+\bar{F}\mathbb{E}x_k$  works from $k=\psi+1$.
Let $\mathbb{V}_k^{F, \bar{F}; \mathbb{V}_\psi^{v_\psi}}$  denote the state
solution for the system $[\mathscr{A}_1^{F, \bar{F}};\mathscr{A}_2^{F},
\mathscr{A}_3^{F}]$ starting from the initial state $v_\psi$. The cost function is redefined as
$$
\hat{\tilde{J}}(F, \bar{F}; \mathbb{V}_\psi^{v_\psi})\triangleq \mathbb{E}  \sum\limits_{k=\psi}^\infty
(\mathbb{V}_k^{F, \bar{F}; \mathbb{V}_\psi^{v_\psi}})'
\left[\begin{array}{ccccc}
\Lambda+\bar{\Lambda}& \bf{0} \\
\bf{0} &\Lambda \end{array}\right]\mathbb{V}_k^{F, \bar{F}; \mathbb{V}_\psi^{v_\psi}}.
$$
In this case, although
$$
\min\limits_{(F, \bar{F})\in \mathcal{F}}\hat{\tilde{J}}(F, \bar{F}; \mathbb{V}_\psi^{v_\psi})\neq \min\limits_{(F, \bar{F})\in \mathcal{F}}{\tilde{J}}(F, \bar{F}; \mathbb{V}_\psi^{z}),
$$
 the optimal solution  ($F^*, \bar{F}^*$) is unchanged.
Besides, for the GLE (\ref{gle1}) in Primal Problem  I,   the terms
$\mathbb{F}_1Z_2\mathbb{F}_1'+
\mathbb{F}_2(Z_1-Z_2)\mathbb{F}_2'$  is replaced by a positive definite
matrix $
\aleph=\sum\limits_{l=1}^r v_lv_l'>0.
$
According to  {Theorem 1 in \cite{chenchen}} and Proposition  \ref{proposition123}, we can derive 
the following  modified Primal Problem $I'$.

{\bf  Problem 5} (Primal  Problem $I'$).
Solve the following  convex minimization problem with variables $\tilde{\mathbbm{S}}_{\mathcal{P}_{I'}}\in\mathscr{S}_{2n+2m}$ and  $(F_{\mathcal{P}_{I'}}\times\bar{F}_{\mathcal{P}_{I'}})\in (\mathcal{R}^{m \times n}\times
\mathcal{R}^{m \times n})$:
\begin{eqnarray}\label{gle3}
\left\{\begin{array}{l}
J_{\mathcal{P}_{I'}}
\! \triangleq \!\\
\inf\limits_{\tilde{\mathbbm{S}}
_{\mathcal{P}_{I'}}\! \in\! \mathscr{S}_{2n+2m}, \! \ F_{\mathcal{P}_{I'}}, \! \bar{F}_{\mathcal{P}_{I'}}\! \in\! \mathcal{R}^{m\times n}}
\! Tr\! \Big(\! \left[\!\begin{array}{ccccc}
\Lambda\! +\! \bar{\Lambda}&\bf{0} \\
\bf{0} &\Lambda \end{array}\!\right]\! \tilde{\mathbbm{S}}_{\mathcal{P}_{I'}}\! \Big),\\
\text{s.t.}\
 \mathscr{A}_1^{F, \bar{F}}\tilde{\mathbbm{S}}_{\mathcal{P}_{I'}} (\mathscr{A}_1^{F, \bar{F}})'
+\mathscr{A}_2^{F}\tilde{\mathbbm{S}}_{\mathcal{P}_{I'}} (\mathscr{A}_2^{F})'
\\
\ \ \ \ \ \ \ \ \ \   +\mathscr{A}_3^{F}\tilde{\mathbbm{S}}_{\mathcal{P}_{I'}}
(\mathscr{A}_3^{F})'+\aleph=\tilde{\mathbbm{S}}_{\mathcal{P}_{I'}},\\
\ \ \ \ \
\tilde{\mathbbm{S}}_{\mathcal{P}_{I'}}>0.
 \end{array}
\right.
\end{eqnarray}
\begin{remark}
Primal Problem $I'$  is equivalent to Primal Problem  I in the sense of
$F^*_{\mathcal{P}_I}=F^*_{\mathcal{P}_{I'}}$, $\bar{F}^*_{\mathcal{P}_I}=\bar{F}^*_{\mathcal{P}_{I'}}$, and $\tilde{\mathbbm{S}}_{\mathcal{P}_{I}}=\tilde{\mathbbm{S}}
_{\mathcal{P}_{I'}}$, if the   initial states $z_{(i)}$ and $v_{(i)}$
  of the state trajectories $x_k$ and $v_k$
  satisfy the following relationship
\begin{align*}\label{inizv}
&z_i=\Big(\left[\begin{array}{cc} A_1\ B_1
\end{array}\right] {v}_i
+\left[\begin{array}{cc} \bar{A}_1\ \bar{B}_1
\end{array}\right]\bar{v}_i\Big)\Big(
\left[\begin{array}{cc} A_1\ B_1
\end{array}\right] {v}_i\nonumber\\
&  \ \ \ \ \ \ \
+\left[\begin{array}{cc} \bar{A}_1\ \bar{B}_1
\end{array}\right]\bar{v}_i\Big)'
+\Big(\left[\begin{array}{cc} A_2\ B_2\end{array}\right]v_i
+\left[\begin{array}{cc} \bar{A}_2\ \bar{B}_2
\end{array}\right]\bar{v}_i\Big)\nonumber\\
& \ \ \ \ \ \ \ \cdot\Big(\left[\begin{array}{cc} A_2\ B_2\end{array}\right]v_i
+\left[\begin{array}{cc} \bar{A}_2\ \bar{B}_2
\end{array}\right]\bar{v}_i\Big)',  \nonumber\\
 & \bar{v}_i=\mathbb{E}{v}_i, \ i=1, 2, \cdots, r.
\end{align*}
\end{remark}
A similar analysis of strong duality   for the Primal Problem  I
 also applies to the  modified Primal Problem $I'$,  with the corresponding dual problem given by
\begin{eqnarray*}
\hat{J}_{\mathcal{D}}:= \sup\limits_{\tilde{\mathfrak{X}}, \tilde{\mathfrak{X}}_0\in\mathscr{S}_{2n+2m}}
\inf\limits_{\tilde{\mathbbm{S}}\in\mathscr{S}_{2n+2m}^+, F, \bar{F}\in\mathcal{F}}
\hat{L}(\tilde{\mathfrak{X}}, \tilde{\mathfrak{X}}_0, F, \bar{F}, \tilde{\mathbbm{S}}),
\end{eqnarray*}
where
\begin{eqnarray*}
&&\hat{L}(\tilde{\mathfrak{X}}, \tilde{\mathfrak{X}}_0, F, \bar{F}, \tilde{\mathbbm{S}})\\
&&=
Tr(\tilde{\mathfrak{X}}\aleph)+
Tr\Big\{
 [(\mathscr{A}_1^{F, \bar{F}})'\tilde{\mathfrak{X}}
 \mathscr{A}_1^{F, \bar{F}}
+(\mathscr{A}_2^{F})'\tilde{\mathfrak{X}}
\mathscr{A}_2^{F}\\
&& \ \ \ \
+(\mathscr{A}_3^{F})'\tilde{\mathfrak{X}}
\mathscr{A}_3^{F}+\left[\begin{array}{ccccc}
\Lambda+\bar{\Lambda}& \bf{0} \\
\bf{0} &\Lambda \end{array}\right]-\tilde{\mathfrak{X}}
 -\tilde{\mathfrak{X}}_0]\tilde{\mathbbm{S}}\Big\}.
\end{eqnarray*}
 This implies that
$J_{\mathcal{P}_{I'}}=\hat{J}_{\mathcal{D}}$.  \\
The following result derives the KKT conditions, which are
both  necessary and sufficient for the optimal solution of the  Problem 5.
\begin{proposition}\label{prpddddd}
For Problem 5, $(\tilde{\mathbbm{S}}^*, F^*, \bar{F}^*)$ is the primary optimal point
and $(\tilde{\mathfrak{X}}^*, \tilde{\mathfrak{X}}_0^*)$ is the dual
optimal point iff $(\tilde{\mathbbm{S}}^*, F^*, \bar{F}^*, \tilde{\mathfrak{X}}^*, \tilde{\mathfrak{X}}_0^*)$ satisfies the following KKT conditions:
\begin{eqnarray}
&& 
 \mathscr{A}_1^{F^*, \bar{F}^*}\tilde{\mathbbm{S}}^*
  (\mathscr{A}_1^{F^*, \bar{F}^*})'
+\mathscr{A}_2^{F^*}\tilde{\mathbbm{S}}^*
 (\mathscr{A}_2^{F^*})'
\nonumber \\
&& +\mathscr{A}_3^{F^*}\tilde{\mathbbm{S}}^*
(\mathscr{A}_3^{F^*})'+\aleph-\tilde{\mathbbm{S}}^*=0,\label{kkt1}\\
&& 
\tilde{\mathbbm{S}}^*>0, \label{kkt2}\\
&& 
 (\mathscr{A}_1^{F^*, \bar{F}^*})'\tilde{\mathfrak{X}}^*
\mathscr{A}_1^{F^*, \bar{F}^*}
+(\mathscr{A}_2^{F^*})'\tilde{\mathfrak{X}}^*
\mathscr{A}_2^{F^*}
\nonumber \\
&&  +(\mathscr{A}_3^{F^*})'\tilde{\mathfrak{X}}^*
\mathscr{A}_3^{F^*}+\left[\begin{array}{ccccc}
\Lambda+\bar{\Lambda}& \bf{0} \\
\bf{0} &\Lambda \end{array}\right]-
\tilde{\mathfrak{X}}^*=0,\label{kkt3} \\
&&  \! 
[\bar{\mathfrak{X}}_{12}^{*'}\! +\! \bar{\mathfrak{X}}_{22}^*\!
(F^*\!+\! \bar{F}^*)]\! \left[\begin{array}{ccc}\! \hat{A}_1\! \! &\! \! \hat{B}_1\! \end{array}
\right]\!\bar{\mathbbm{S}}^* \! \left[\begin{array}{ccc}\!\hat{A}_1\! \! &\! \! \hat{B}_1\!\!  \end{array}
\right]'\!\! =\!\!  0,\label{kkt4} \\
 && ({\mathfrak{X}}_{12}^{*'}+{\mathfrak{X}}_{22}^*
F^*)\Psi=0  \label{kkt5}
\end{eqnarray}
with
\begin{eqnarray*}
&&\Psi \!=\left[\begin{array}{ccc}\! {A}_1\!\!&\!\! {B}_1\!\end{array}
\right] {\mathbbm{S}}^* \!\left[\begin{array}{ccc} \!{A}_1\!\!& \!\! {B}_1\!\end{array}
\right]'
\!+\!
\left[\begin{array}{ccc}\! {A}_2\!\!&\!\! {B}_2\end{array}
\right] \!{\mathbbm{S}}^*\! \left[\begin{array}{ccc}\! {A}_2\!\!&\!\! {B}_2\!\end{array}
\right]'\\
&& \ \ \ \ \ \ \ \ \
+\left[\begin{array}{ccc}\!\hat{A}_2\!\!&\!\!\hat{B}_2\!\end{array}
\right]\!\bar{\mathbbm{S}}^* \! \left[\begin{array}{ccc}\!\hat{A}_2\!\!&\!\!\hat{B}_2\!\end{array}
\right]'.
\end{eqnarray*}
\end{proposition}
{\textit{Proof:}} The proof is provided in Appendix F.\\
Based on the KKT conditions (\ref{kkt1})-(\ref{kkt5}),
 the primary optimal point  ($F, \bar{F}$) and
 dual optimal point  $\tilde{\mathfrak{X}}$
  can be obtained by iteratively solving
 the stationary conditions (\ref{kkt3})-(\ref{kkt5}), which is proposed in
 Algorithm 2. Moreover, the convergence  analysis for the PD algorithm
  is given in the following result.
  \begin{theorem}\label{theorem41}
  In Algorithm 2, the sequences
  $\{\tilde{\mathfrak{X}}_i\}_{i=0}^\infty$ and $\{F^i, \bar{F}^i\}_{i=0}^\infty$ satisfy
  $$
  \lim\limits_{i\rightarrow \infty}
  \tilde{\mathfrak{X}}_i =\left[\begin{array}{ccc}
  \bar{\mathfrak{X}}^* & \bf{0} \\ \bf{0} &  {\mathfrak{X}}^*
  \end{array}\right], \
   \lim\limits_{i\rightarrow \infty}
   {F}^i =F^*, \
   \lim\limits_{i\rightarrow \infty}
   {\bar{F}}^i =\bar{F}^*,
  $$
  where
 $ {\mathfrak{X}}^*$  and  $\bar{\mathfrak{X}}^*$ are
  defined in  (\ref{qfunction1}) and (\ref{qfunction2}), respectively.
  $F^*$ and $\bar{F}^*$ are optimal feedback stabilizing control
gains as in (\ref{gain1})-(\ref{gain2}).
 \end{theorem}
 {\textit{Proof:}} The proof is provided in Appendix G.\\
  \begin{algorithm}
  \caption{Model Based PD Algorithm}
  \begin{algorithmic}[1]  
          \State {\bf Initialization:} Select any stabilizer $(F^0, \bar{F}^0)\in (\mathcal{F}\times \mathcal{F})$. Let   the convergence tolerance
          $\varepsilon >0$  and   $i=0$;
           \State  {\bf Repeat Steps $3-6$;}
            \State  {\bf Dual Update:} Solve $\tilde{\mathfrak{X}}_i\triangleq
            \left[\begin{array}{ccccc}
\bar{\mathfrak{X}}_i& \bf{0} \\
 \bf{0} &{\mathfrak{X}}_i \end{array}\right]
            $ from
            the   equation:
            \begin{eqnarray}
            &(\mathscr{A}_1^{F^i, \bar{F}^i})'\tilde{\mathfrak{X}}_i
\mathscr{A}_1^{F^i, \bar{F}^i}
+(\mathscr{A}_2^{F^i, \bar{F}^i})'\tilde{\mathfrak{X}}_i
\mathscr{A}_2^{F^i, \bar{F}^i}\nonumber \\
&
+(\mathscr{A}_3^{F^i, \bar{F}^i})'\tilde{\mathfrak{X}}_i
\mathscr{A}_3^{F^i, \bar{F}^i}  +\left[\begin{array}{ccccc}
\Lambda+\bar{\Lambda}& \bf{0} \\
 \bf{0} &\Lambda \end{array}\right]=
\tilde{\mathfrak{X}}_i; \label{dualupdate}
              \end{eqnarray}
                         \State  {\bf Primal Update:} Update control gains from
               \begin{equation}\label{dualupdate11}
             F^{i+1}=
     ({\mathfrak{X}}^{22*}_i)^{-1}
     ({\mathfrak{X}}^{12*}_i)',
             \end{equation}
             and
              \begin{equation}\label{dualupdate22}
              \bar{F}^{i+1}=
              (\bar{\mathfrak{X}}_i^{22*})^{-1}
     (\bar{\mathfrak{X}}_i^{12*})'-
      ({\mathfrak{X}}_i^{22*})^{-1}
     ({\mathfrak{X}}_i^{12*})' ;
             \end{equation}
                           \State  $i\leftarrow i+1$;
             \State {\bf Until} $\|F^i-F^{i+1}\|+\|\bar{F}^i-\bar{F}^{i+1}\|\leq \varepsilon$.
    \end{algorithmic}
\end{algorithm}
Defining a new matrix $\tilde{\mathbbm{S}}^M(F, \bar{F})$
by truncating the time horizon of $\tilde{{\mathbbm{S}}}$ at finite
$M>0$:
$$
\tilde{\mathbbm{S}}^M(F, \bar{F})\triangleq \sum_{l=1}^r\sum_{k=\psi}^M
\mathbb{E}[(\mathbb{V}_k^{F, \bar{F}; \mathbb{V}_0^{v_l}})
(\mathbb{V}_k^{F, \bar{F}; \mathbb{V}_0^{v_l}})'].
$$
Obviously, $\tilde{\mathbbm{S}}^M(F, \bar{F})>0$ holds for the positive definiteness of
initial state matrix $\aleph=\sum\limits_{l=1}^r v_lv_l'>0.$  For the KKT condtion (\ref{kkt3}), a simple calculation leads to
\begin{eqnarray}
&&
\mathbbm{W}^M(F, \bar{F})\tilde{\mathfrak{X}}
(\mathbbm{W}^M(F, \bar{F}))'
+\tilde{\mathbbm{S}}^M(F, \bar{F})(\mathscr{A}_2^{F})'\tilde{\mathfrak{X}}
\mathscr{A}_2^{F}\nonumber\\
&&\cdot \tilde{\mathbbm{S}}^M(F, \bar{F})
+\tilde{\mathbbm{S}}^M(F, \bar{F})(\mathscr{A}_3^{F})'\tilde{\mathfrak{X}}
\mathscr{A}_3^{F}\tilde{\mathbbm{S}}^M  (F, \bar{F}) \nonumber\\
&&+\tilde{\mathbbm{S}}^M(F, \bar{F})  \Big(\left[\begin{array}{ccccc}
\Lambda+\bar{\Lambda}&\bf{0} \\
\bf{0} &\Lambda \end{array}\right]-\tilde{\mathfrak{X}}\Big)\tilde{\mathbbm{S}}^M(F, \bar{F})=0  \label{kkt33}
\end{eqnarray}
with $\mathbbm{W}^M(F, \bar{F})\triangleq \sum_{l=1}^r\sum_{k=\psi}^M
\mathbb{E}[(\mathbb{V}_k^{F, \bar{F}; \mathbb{V}_\psi^{v_l}})
(\mathbb{V}_{k+1}^{F, \bar{F}; \mathbb{V}_\psi^{v_l}})']$.
Based on the modified KKT condition (\ref{kkt33}), the dual update
$\tilde{\mathfrak{X}}_i$ can be solved without the knowledge of system matrices
$\mathscr{A}_1^{F, \bar{F}}$, i.e., the system drift term parameters $A_1, \bar{A}_1, B_1$ and $\bar{B}_1$. The partial model-free algorithm has been given as in Algorithm 3.
  \begin{algorithm}
  \caption{Partially Model-Free  PD Algorithm}
  \begin{algorithmic}[1]  
          \State {\bf Initialization:} Select any stabilizer $(F^0, \bar{F}^0)\in
           (\mathcal{F}\times \mathcal{F})$. Let   the convergence tolerance
          $\varepsilon >0$  and   $i=0$;
           \State  {\bf Repeat Steps $3-7$;}
            \State  {\bf  }
            By calculating the mean value $\mathbb{E}[x_k(F^i, \bar{F}^i)]\thickapprox \frac{1}{H}\sum_{h=1}^H x_{k, h}(F^i, \bar{F}^i)$,
            solve $\tilde{\mathbbm{S}}^M(F^i, \bar{F}^i)$ and $\mathbbm{W}^M(F^i, \bar{F}^i)$;
              \State  {\bf Dual Update:} Calculate $\tilde{\mathfrak{X}}^i$ by
              solving (\ref{kkt33}) with the iteration $F^i$ and $\bar{F}^i$;
                         \State  {\bf Primal Update:} Update control gains from
               \begin{equation*}
             F^{i+1}=
     ({\mathfrak{X}}^{22*}_i)^{-1}
     ({\mathfrak{X}}^{12*}_i)',
             \end{equation*}
             and
              \begin{equation*}
              \bar{F}^{i+1}=
              (\bar{\mathfrak{X}}_i^{22*})^{-1}
     (\bar{\mathfrak{X}}_i^{12*})'-
      ({\mathfrak{X}}_i^{22*})^{-1}
     ({\mathfrak{X}}_i^{12*})' ;
             \end{equation*}
                           \State  $i\leftarrow i+1$;
             \State {\bf Until} $\|F^i-F^{i+1}\|+\|\bar{F}^i-\bar{F}^{i+1}\|\leq \varepsilon$.
    \end{algorithmic}
\end{algorithm}

\section{Simulation}
The solvability of the MF-SLQR problem is   challenging due to the coupled AREs and the presence of MF terms $\mathbb{E}x_k$ and $\mathbb{E}u_k$ in the stochastic system. Existing studies on the MF-SLQR problem  have primarily focused on finite horizon settings \cite{elliott2013} or 1-dimensional dynamic systems \cite{chenchen,qqy}.  In this section, we demonstrate the effectiveness of the model-based and partially model-free PD-PL algorithms through simulations on a 3-dimensional system \eqref{sys1}, where the system matrices are selected from \cite{elliott2013} and are shown as follows:
 \begin{figure*}[htbp]
\centering
\subfigure[\ ]{
    \includegraphics[width=0.3\textwidth]{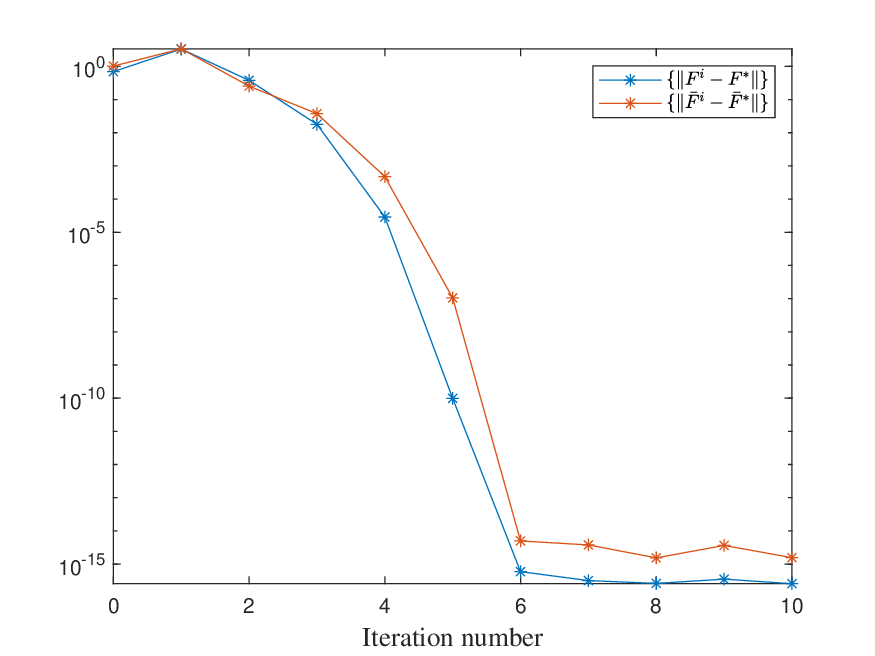}
}
\subfigure[ ]{
    \includegraphics[width=0.3\textwidth]{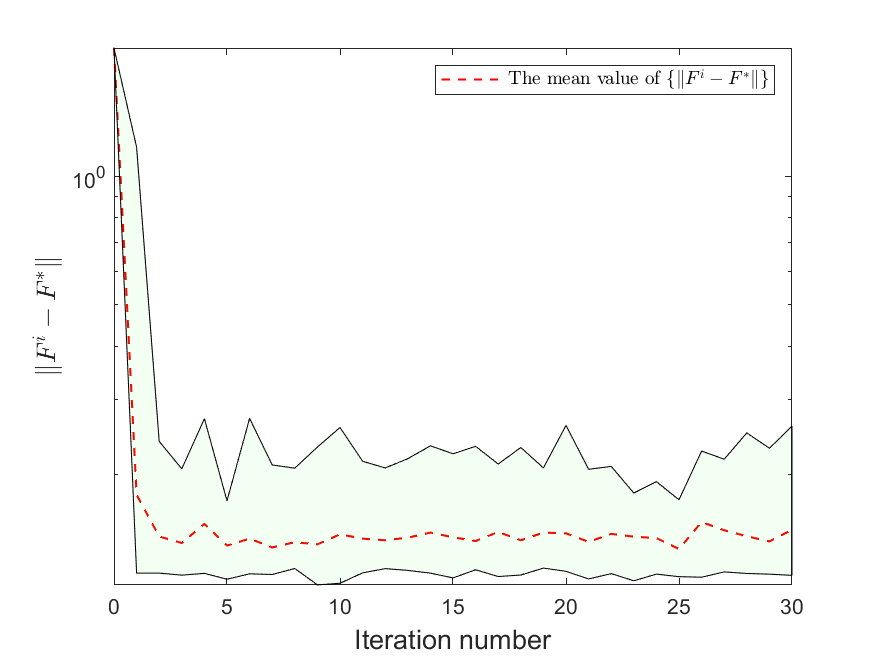}
}
\subfigure[ ]{
    \includegraphics[width=0.3\textwidth]{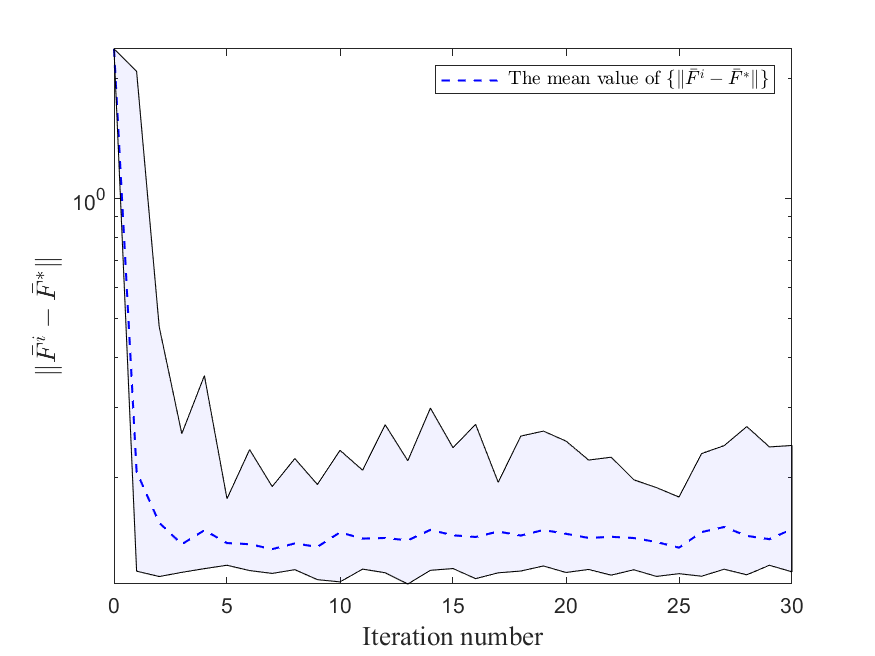}
}
\subfigure[ ]{
    \includegraphics[width=0.3\textwidth]{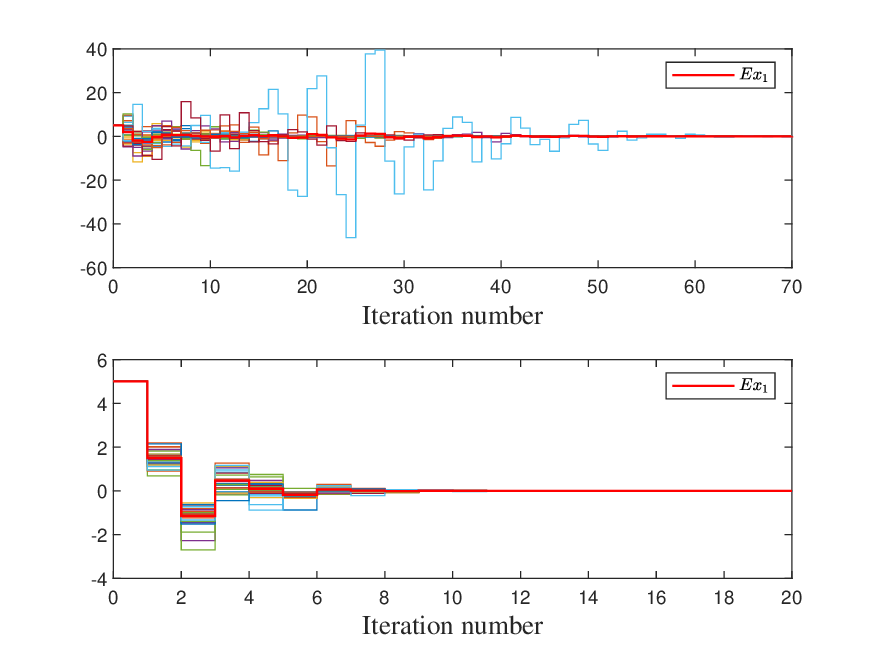}
}
\subfigure[ ]{
    \includegraphics[width=0.3\textwidth]{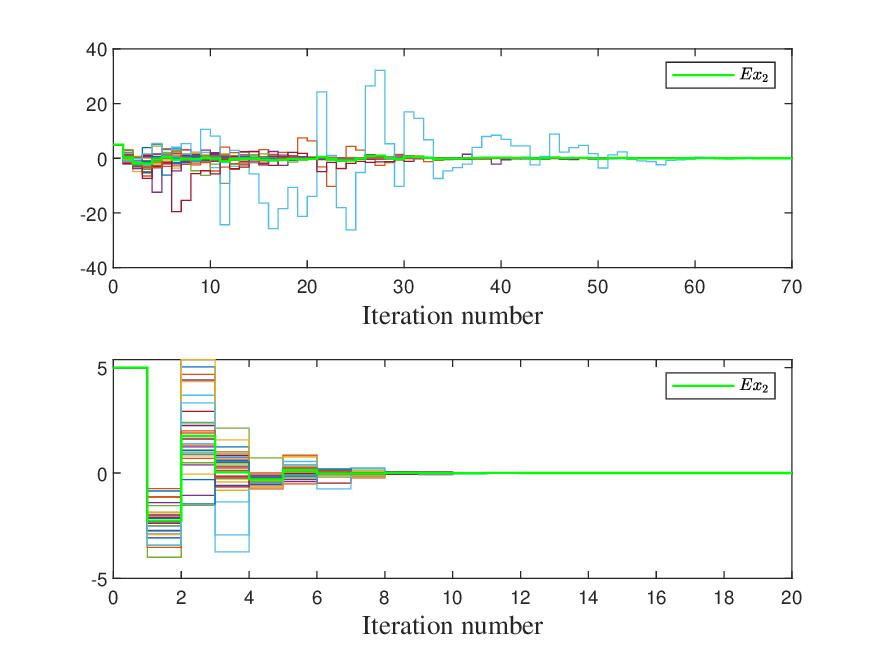}
}
\subfigure[ ]{
    \includegraphics[width=0.3\textwidth]{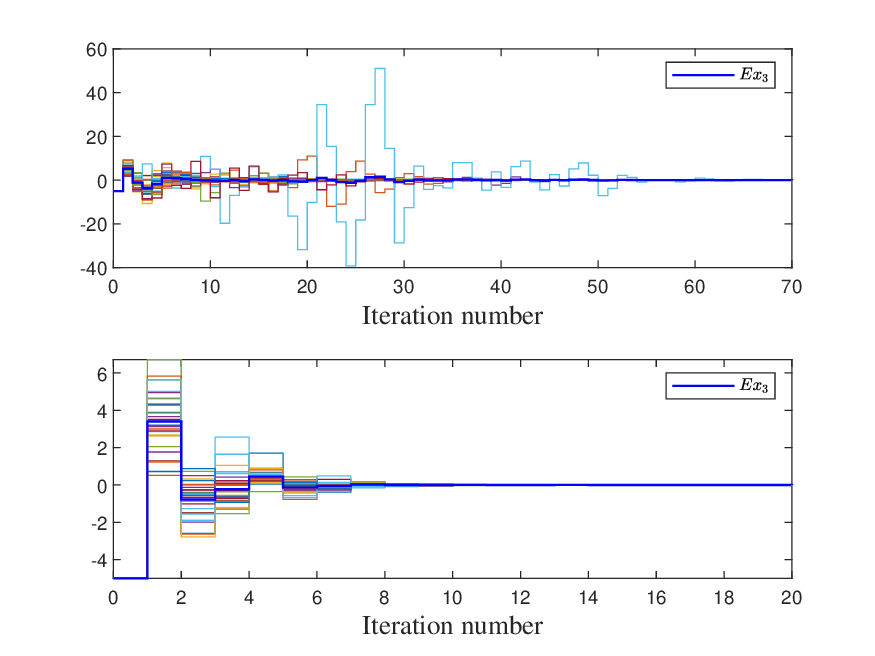}
}
\caption{(a): Simulation results for model-based PD Algorithm; (b-c): Simulation results for partially model-free PD algorithm; \ \ \ \ \ \ (d-f): Measured data from simulations.}
\label{fig:simulation_results}
\end{figure*}
\begin{eqnarray*}
&&A_1\!=\!\begin{bmatrix}0.2&0.4&0.2\\0&0.2&0.6\\0.6&0.4&0.2\end{bmatrix},\!
\bar{A}_1\!=\!\begin{bmatrix}0.3&0.4&0.2\\0&0.2&0.7\\0.6&0.5&0.2
\end{bmatrix}, \!
\\
&&C_1\!=\!\begin{bmatrix}0.2&0.4&0.6\\0.4&0.2&0.6\\0.2&0.4&0.2\end{bmatrix},\!
\bar{C}_1\!=\!\begin{bmatrix}0.3&0.4&0.6\\0.4&0.3&0.6\\0.2&0.4&0.3
\end{bmatrix},\!
\\
&&B_1\!=\!\begin{bmatrix}0.4&0.2\\0.2&0.4\\0.3&0.3\end{bmatrix},
\bar{B}_1\!=\!\begin{bmatrix}0.5&0.2\\0.2&0.5\\0.2&0.3\end{bmatrix},\\
&&D_1=\begin{bmatrix}0.2&0.6\\0.6&0.4\\0.3&0.1\end{bmatrix},
\bar{D}_1=\begin{bmatrix}0.3&0.5\\0.5&0.4\\0.3&0.3\end{bmatrix},\\
&&Q=diag(0,1.5,1),\ \bar{Q}=diag(1,1,0),\\
&&R=diag(1,1)\ \bar{R}=diag(1.5,1).
\end{eqnarray*}
\begin{table}[htbp]
\centering
\caption{Values of $F^0$, $\bar{F}^0$, $F^*$ and $\bar{F}^*$}
\begin{tabular}{|c|c|c|c|}
\hline
\textbf{Parameter} & \textbf{Value} \\ \hline
$F^0$ & $\left[\begin{array}{ccc}-0.73&-1.50&-1.30\\ 0.23&0.97&-0.26\end{array}\right]$ \\ \hline
$F^*$ & $\left[\begin{array}{ccc}-0.34  & -0.32  & -0.42\\ -0.24  & -0.29 &  -0.49\end{array}\right]$ \\ \hline
$\bar{F}^0$ & $\left[\begin{array}{ccc}1.21&1.81&0.66\\ -0.18&-0.16&-1.50\end{array}\right]$ \\ \hline
$\bar{F}^*$ & $\left[\begin{array}{ccc}-0.32 &-0.42  & -0.34\\ -0.31  & -0.43 &  -0.77\end{array}\right]$ \\ \hline
\end{tabular}
\label{tab:final_values}
\end{table}
Select  the initial stabilizing control gains  $F^0$ and $\bar{F}^0$ by observing whether they can stabilize a diverging state. The analytical solutions for the optimal feedback stabilizing control gains $F^*$ and $\bar{F}^*$  are obtained by the PI  algorithm, i.e., Algorithm 1.   The values are given in Table 1.
For comparison, we apply Algorithms 2 and 3 to a system based on \cite{elliott2013}. Unlike in \cite{elliott2013}, the current system is only partially known, with the drift terms $A_1$, $\bar{A}_1$, $B_1$, and $\bar{B}_1$  unspecified.  The time horizon is extended from horizon $k=3$ to  $k=\infty$.
In order to  satisfy the initial state condition setting in Problem 1, we randomly choose the initial state given as follows:
{\small $\begin{array}{llllllllllllllll}
&\ \ \ \ \left\{Ex^i_0 \mid i = 1, \cdots, 20\right\}\\
&=\left\{
\begin{bmatrix}9.59 \\1.55 \\3.15\end{bmatrix},
\begin{bmatrix}8.74 \\2.85 \\2.85\end{bmatrix},
\begin{bmatrix}3.71 \\1.71 \\7.35\end{bmatrix},
\begin{bmatrix}6.90 \\2.52 \\5.17\end{bmatrix},
\begin{bmatrix}2.54 \\2.49 \\6.13\end{bmatrix},
\begin{bmatrix}5.04 \\1.32 \\1.67\end{bmatrix},\right.\\ &\ \ \ \ \ \ \left.
\begin{bmatrix}1.32 \\7.40 \\6.21\end{bmatrix},
\begin{bmatrix}7.58 \\7.55 \\0.93\end{bmatrix},
\begin{bmatrix}2.55 \\2.97 \\2.34\end{bmatrix},
\begin{bmatrix}1.74 \\1.72 \\0.72\end{bmatrix},
\begin{bmatrix}1.45 \\5.21 \\0.30\end{bmatrix},
\begin{bmatrix}4.90 \\0.70 \\2.58\end{bmatrix},\right.
\end{array}$}
{\small $\begin{array}{llllllllllllllll}
&\ \ \ \ \ \ 
\left.
\begin{bmatrix}5.65 \\6.15 \\3.05\end{bmatrix},
\begin{bmatrix}6.02 \\5.13 \\7.47\end{bmatrix},
\begin{bmatrix}5.31 \\6.65 \\6.87\end{bmatrix},
\begin{bmatrix}3.62 \\0.09 \\1.45\end{bmatrix},
\begin{bmatrix}8.98 \\3.18 \\9.92\end{bmatrix},
\begin{bmatrix}2.85 \\3.58 \\2.57\end{bmatrix},\right.\\ &\ \ \ \ \ \ \left.
\begin{bmatrix}7.35 \\3.72 \\7.35\end{bmatrix},
\begin{bmatrix}9.50 \\4.54 \\4.02\end{bmatrix}
\right\}
\end{array}$}
 and
{\small $\begin{array}{llllllllllllllll}
&\ \ \ \ \left\{x^i_0-Ex^i_0 \mid i = 1, \cdots, 20\right\}\\
&=\left\{
\begin{bmatrix}-5.50 \\6.65 \\11.47\end{bmatrix},
\begin{bmatrix}-3.05 \\3.85 \\-3.65\end{bmatrix},
\begin{bmatrix}8.58 \\12.14 \\7.27\end{bmatrix},
\begin{bmatrix}-4.92 \\7.72 \\-3.98\end{bmatrix},
\begin{bmatrix}11.51 \\4.60 \\-2.93\end{bmatrix},\right.
\\
   &\ \ \ \ \ \ \left.
   \begin{bmatrix}8.66 \\1.77 \\-1.64\end{bmatrix},
\begin{bmatrix}3.69 \\6.37 \\0.35\end{bmatrix},
\begin{bmatrix}-3.62 \\-4.21 \\0.29\end{bmatrix},
\begin{bmatrix}-1.04 \\10.41 \\6.30\end{bmatrix},
\begin{bmatrix}3.57 \\3.39 \\7.02\end{bmatrix},\right.\\
   &\ \ \ \ \ \ \left.
   \begin{bmatrix}8.54 \\-2.21 \\3.30\end{bmatrix},
\begin{bmatrix}-0.64 \\1.39 \\3.37\end{bmatrix},
\begin{bmatrix}7.51 \\4.67 \\7.87\end{bmatrix},
\begin{bmatrix}5.09 \\-1.92 \\0.92\end{bmatrix},
\begin{bmatrix}1.49 \\-4.78 \\2.09\end{bmatrix},\right.
\end{array}$}
{\small $\begin{array}{lllllllllllllllll}
 &\ \ \ \ \ \ \left.
\begin{bmatrix}11.17 \\11.02\\11.89\end{bmatrix},
\begin{bmatrix}5.00 \\11.70 \\-8.55\end{bmatrix},
\begin{bmatrix}5.80 \\1.03 \\10.92\end{bmatrix},
\begin{bmatrix}-3.16 \\0.33 \\1.54\end{bmatrix},
\begin{bmatrix}-6.89 \\9.02 \\6.52\end{bmatrix}.
\right\}
\end{array}$}

By implementing the model-based PD Algorithm, Fig. 1-(a) demonstrates that the control gains converge to their optimal values with an accuracy of $10^{-15}$   at  iteration $ i=10$.  Using the partially model-free PD algorithm, the convergence process terminates after $i=30$ iterations. The learning trajectories of $F^i$ and $\bar{F}^i$ during this process are illustrated in Figs. 1-(b) and (c), respectively. Here, the Monto Carlo method  is adopted by experimenting 30 times. The state trajectories generated by the control gains before and after the learning process are simulated in Figs. 1-(d-f), where in Figs. 1-(d-f), the upper subfigures represent the state trajectories collected using the initial control, while the lower subfigures represent the state collected using the control obtained from Algorithm 3. The results indicate that the learned optimal trajectories converge to zero significantly faster.

We compare our proposed method with a model-based approach consisting of two steps:
(i) Estimate the system parameters $({\mathcal{A}_1},{\mathcal{B}}_1,{\bar{\mathcal{A}}_1},{\bar{\mathcal{B}}_1})$  
 from trajectory data   $\{\mathbb{E}x_k\}_{k=0}^{r}$ and $\{x_k\}_{k=0}^{r}$ (  $r>n+m$)  using least squares; (ii)
Solve GAREs \eqref{gare1}--\eqref{gare2} with these estimates via the SDP method. The estimation procedure is as shown in \cite{lina2022}. Using this procedure, the feedback gain computed via \eqref{gain2} with estimated parameters satisfies
$\|\hat{\mathbb{F}}-\hat{F}^*\|=3.751\times 10^{-1}$, while the proposed learning-based method achieves $\|\hat{F}^{5}-\hat{F}^*\|=1.013\times 10^{-1}$.  Thus, the partially model-free approach is significantly more accurate than the model-based method.

\section{Conclusion}
This article has adopted the primal-dual method to research  the PL algorithm for designing the optimal controller for the discrete-time stochastic LQ systems with the MF terms.
The  primal and dual optimization problems have been constructed to reformulate the dynamic MF-SLQR problem into the new nonconvex optimization problems, which can be researched by the  semidefinite programming technique, and strong duality theorem is established.  The uniform  convergence of the PL algorithm has been theoretically analyzed, including the fact that all the learned policies are still   stabilizing. In our future research, we will generalize  the  PD-PL  framework to solve  the risk-constrained MF-SLQR problem and explore other forms of stability, such as input-to-state stability,  to further enhance the algorithm's  robustness.

 \section*{Appendix A}
 {\textit{Proof of Theorem \ref{th22211}.}}\\
$(1)$.  To begin, we will establish the following   iterative formula:
 \begin{eqnarray}\label{pe21}
         &&\left[\begin{array}{cc}
           \bar{P}^i&  \bf{0} \\ \bf{0}   &  P^i\end{array}\right] \nonumber\\
           &&  =  (\mathcal{A}_1^{F^{i+1}, \bar{F}^{i+1}})'\!
              \left[\begin{array}{cc}
             \! \bar{P}^i  \!\! &\!\! \bf{0}\! \\ \!\bf{0}\!\!  &\!\!  P^i\!\end{array}\right](\mathcal{A}_1^{F^{i+1}, \bar{F}^{i+1}})\!+\!\widetilde{\tilde{\mathcal{Q}}}^{i+1}
           \!+\!(\mathcal{A}_2^{F^{i+1}})' \nonumber\\
           &&
                 \cdot  \! \left[\! \begin{array}{cc}
             \!\bar{P}^i\!\! &\!\! \bf{0} \!\\ \!\bf{0} \!\! &\!\!P^i\!\end{array}\! \right]\!
             (\! \mathcal{A}_2^{F^{i+1}}\! )
            \!  + \! (\mathcal{A}_3^{F^{i+1}, \bar{F}^{i+1}}\! )' \! \left[\! \begin{array}{cc}\!
             \bar{P}^i  \!\!& \!\!\bf{0} \!\\ \!\bf{0}\!\!  &\!\!  P^i\!\end{array}\! \right]\!
             (\!\mathcal{A}_3^{F^{i\!+\!1},  \bar{F}^{i\!+\!1}}\! )
                        \end{eqnarray}
           with
           \begin{eqnarray*}
         &&\widetilde{\tilde{\mathcal{Q}}}^{F^{i+1}, \bar{F}^{i+1}}\\
           &&=\tilde{\mathcal{Q}}^{i+1}  \! +\! \left[\begin{array}{cc}
             \!\Delta \hat{F}^{i} \!\! &\!\! \bf{0}\! \\ \!\bf{0}  \!\! & \!\! \Delta F^{i}\!\end{array}\right]'\!
                         \left[\begin{array}{cc}
            \!\Upsilon^{P^i, \bar{P}^i} \!\!  &\!\! \bf{0}\! \\\! \bf{0} \!\!& \!\! \Upsilon^{P^i}\!\end{array}\right]\!
            \left[\begin{array}{cc}\!
             \Delta \hat{F}^{i}\!\!  & \!\!\bf{0} \!\\ \!\bf{0} \!\! &\!\! \Delta F^{i}\end{array}\right].
            \end{eqnarray*}
Denote
\begin{eqnarray*}
&&\Gamma_{\mathcal{A}_1^{F^i, \bar{F}^i},
\mathcal{A}_2^{F^i}, \mathcal{A}_3^{F^i, \bar{F}^i}}\\
&& :=(\mathcal{A}_1^{F^i, \bar{F}^i})'\!
              \left[\begin{array}{cc}
           \!  \bar{P}^i  \!\!& \!\!\bf{0}\! \\ \!\bf{0}\!\!  &  \!\! P^i\!\end{array}\right]\!(\mathcal{A}_1^{F^i, \bar{F}^i})
         \! +\!(\mathcal{A}_2^{F^i})' \!\left[\begin{array}{cc}
            \! \bar{P}^i \!\!&\!\! \bf{0} \!\\ \!\bf{0}\!\!&\!\!P^i\!\end{array}\right]\! (\mathcal{A}_2^{F^i}) \!\\
            && \ \ \ \
            +(\mathcal{A}_3^{F^i, \bar{F}^i})' \!\left[\begin{array}{cc}
             \!\bar{P}^i  \!\!&\!\! \bf{0}\! \\ \!\bf{0}\!\!    &  \!\! P^i\!\end{array}\right]\!
             (\mathcal{A}_3^{F^i, \bar{F}^i}),\\
             &&
             \Gamma_{\mathcal{A}_1, \mathcal{A}_2, \mathcal{A}_3}\\
             &&:=
             \mathcal{A}_1'\!\left[\begin{array}{cc}
             \!\bar{P}^i\!\! &\!\!  \bf{0}\! \\\! \bf{0} \!\! &\!\!P^i\!\end{array}\right]\!
             \left[\begin{array}{cc}
             \!\hat{B}_1\!\! & \!\! \bf{0}\! \\\! \bf{0} \!\!&\!\!B_1\!\end{array}\right]\!
             +\!\mathcal{A}_2'\!\left[\begin{array}{cc}
             \!\bar{P}^i\!\! &\!\!  \bf{0} \!\\\! \bf{0} \!\! &\!\!P^i\!\end{array}\right]\!
             \left[\begin{array}{cc}
            \! \bf{0}\!\!  & \!\!  \bf{0}\! \\ \!\bf{0}   \!\!&\!\!B_2\!\end{array}\right] \\
              && \ \ \ \
             +\mathcal{A}_3'
             \left[\begin{array}{cc}
             \bar{P}^i &  \bf{0} \\ \bf{0}  &P^i\end{array}\right]
             \left[\begin{array}{cc}
            \bf{0}  &  \bf{0} \\ \hat{B}_2  & \bf{0} \end{array}\right],\\
            &&\Gamma_{B_1, \hat{B}_1, B_2, \hat{B}_2}\\
&&:=
             \left[\begin{array}{cc}
             \!\hat{B}_1\!\! &\!\! \bf{0} \!\\\! \bf{0}  \!\!&\!\!B_1\!\end{array}\right]'\!
             \left[\begin{array}{cc}
             \!\bar{P}^i \!\!&\!\!  \bf{0}\! \\ \!\bf{0}\!\! &\!\!P^i\!\end{array}\right]\!
             \left[\begin{array}{cc}
            \! \hat{B}_1 \!\!&  \!\!\bf{0} \!\\\! \bf{0}\!\!  &\!\!B_1\!\end{array}\right]\!
             +\! \left[\begin{array}{cc}
            \! \bf{0}\!\! &\!\!  \bf{0}\! \\ \!\bf{0} \!\!&\!\!B_2\!\end{array}\right]'
            \\
              && \ \ \ \
            \! \cdot \left[\begin{array}{cc}\!
             \bar{P}^i \!\!&\!\!  \bf{0} \!\\ \!\bf{0} \!\! &\!\!P^i\!\end{array}\right]\! \left[\begin{array}{cc}
           \! \bf{0}\!\! &\!\! \bf{0}\! \\ \!\bf{0}  \!\!&\!\!B_2\!\end{array}\right]\! +\!
              \left[\begin{array}{cc}
           \!  \bf{0}\!\!  & \!\! \bf{0}\! \\ \!\hat{B}_2 \!\! &\!\! \bf{0} \! \end{array}\right]'  \!\left[\begin{array}{cc}
           \!  \bar{P}^i \!\!& \!\! \bf{0} \!\\ \!\bf{0} \!\!  &\!\!P^i\!\end{array}\right]\!
             \left[\begin{array}{cc}
            \!\bf{0} \!\! &\!\!  \bf{0}\! \\\! \hat{B}_2  \!\!&\!\!\bf{0} \! \end{array}\right].
              \end{eqnarray*}
  Considering policy update equality (\ref{pe2}), there is
\begin{eqnarray*}
&&\Gamma_{\mathcal{A}_1^{F^{i+1}, \bar{F}^{i+1}},
\mathcal{A}_2^{F^{i+1}}, \mathcal{A}_3^{F^{i+1}, \bar{F}^{i+1}}}\!- \! \left[\begin{array}{cc}
           \!  \bar{P}^i\!\!  & \!\!\bf{0} \!\\\! \bf{0}  \!\!&\!\!  P^i\!\end{array}\right]\\
             &&
             =\Gamma_{\mathcal{A}_1^{F^{i+1}, \bar{F}^{i+1}},
\mathcal{A}_2^{F^{i+1}}, \mathcal{A}_3^{F^{i+1}, \bar{F}^{i+1}}}\!-\!
           \Gamma_{\mathcal{A}_1^{F^{i}, \bar{F}^{i}},
\mathcal{A}_2^{F^{i}}, \mathcal{A}_3^{F^{i}, \bar{F}^{i}}}\\
&&  - \tilde{\mathcal{Q}}^{F^{i}, \bar{F}^{i}}\\
             &&
             =\!-\Gamma_{\mathcal{A}_1, \mathcal{A}_2, \mathcal{A}_3}\left[\begin{array}{cc}
           \!  \Delta \hat{F}^{i}\!\!  &\!\!\bf{0}\! \\ \! \bf{0} \!\! &\!\! \Delta F^{i}\!\end{array}\right]\!
             -\!
             \left[\begin{array}{cc}
             \!\Delta \hat{F}^{i} \!\! &\!\! \bf{0} \!\\ \!\bf{0} \!\!  &\!\! \Delta F^{i}\!\end{array}\right]'\!
             \Gamma_{\mathcal{A}_1, \mathcal{A}_2, \mathcal{A}_3}'\\
              &&  \!- \!
             \left[\begin{array}{cc}
              \!\hat{Q}  \!\! &\!\!   \bf{0}  \!\\  \!\bf{0}  \!\!   & \!\! Q \!\end{array}\right]
             \! + \!
              \left[\begin{array}{cc}
            \! \hat{F}^{i+1}  \!\! & \!\!\bf{0} \! \\ \! \bf{0} \! \!  & \!  \! F^{i+1} \!\end{array}\right]'\!
           \Gamma_{B_1, \hat{B}_1, B_2, \hat{B}_2}\!
           \left[\begin{array}{cc}
           \! \hat{F}^{i+1} \!\! &\!\! \bf{0}\! \\ \!\bf{0} \!\! & \!\!  F^{i+1}\!\end{array}\right]
           \end{eqnarray*}
           \begin{eqnarray*}
             && 
          \!  -
            \! \left[\begin{array}{cc}
          \!  \hat{F}^{i} \!\! & \!\!\bf{0} \!\\ \!\bf{0}  \!\! &\!\!   F^{i}\!\end{array}\right]'\!
            \Big(\!
            \left[\!\begin{array}{cc}
            \!\hat{R}\!\! &\!\! \bf{0}\! \\\! \bf{0} \!\! &\!\!   R \! \end{array}\right]\!+\!
             \Gamma_{B_1, \hat{B}_1, B_2, \hat{B}_2}
           \! \Big)\!\left[\begin{array}{cc}
           \! \hat{F}^{i}  \!\!&\!\! \bf{0}\! \\\! \bf{0}\!  \!&\!\!   F^{i}\!\end{array}\right].
           \end{eqnarray*}
Using policy improvements (\ref{pim1})-(\ref{pim2}), we have
\begin{eqnarray*}
&&-\Gamma_{\mathcal{A}_1, \mathcal{A}_2, \mathcal{A}_3}\\
&&
             =
\left[\begin{array}{cc}
           \! \hat{F}^{i+1} \!\! & \!\!\bf{0} \!\\\! \bf{0} \!\! &\!\!   F^{i+1}\!\end{array}\right]'\!
             \Big(\!
            \left[\begin{array}{cc}
           \! \hat{R} \!\!&  \!\!\bf{0}\!\! \\ \!\bf{0} \!\! &\!\!   R \! \end{array}\right]\!+\!
             \Gamma_{B_1, \hat{B}_1, B_2, \hat{B}_2}\!
            \Big)\\
            &&=\left[\begin{array}{cc}
           \! \hat{F}^{i+1}  \!\!& \!\! \bf{0} \!\\\! \bf{0} \!\! &  \!\! F^{i+1}\!\end{array}\right]'\!
            \left[\begin{array}{cc}
           \! \Upsilon^{P^i, \bar{P}^i}  \!\!&\!\! \bf{0}\! \\ \!\bf{0} \!\! &\!\!  \Upsilon^{P^i}\!\end{array}\right].
            \end{eqnarray*}
After simple calculation, the iterative equation (\ref{pe21}) is obtained.
 Based on the recursions (\ref{pe2}) and  (\ref{pe21}), we know that
 $(F^i, \bar{F}^i)$ are the feedback stabilizing gain matrices, then so are
 $(F^{i+1}, \bar{F}^{i+1})$.  \\
$(2)$.
Based on   {\cite{zwh2004}}, 
if
  $[A_1, \bar{A}_1; A_2, \bar{A}_2|{\mathcal{Q}}]$ is exactly observable,
  then so is
$
[\mathcal{A}_1^{F^i, \bar{F}^i},
\mathcal{A}_2^{F^i}, \mathcal{A}_3^{F^i, \bar{F}^i}|\tilde{\mathcal{Q}}^{F^i, \bar{F}^i}]
$.
Since ${(F^i, \bar{F}^i)}_{i=1}^\infty$ forms a stabilizing gain sequence,
according to Theorem 3 in \cite{qqy},
 the policy evaluation equation (\ref{pe2})
 has a unique solution  $\left[\begin{array}{cc}
           \bar{P}^i\!\!& \!\!\bf{0} \\ \bf{0} \!\!&\!\!P^i\end{array}\right]
           \in \mathscr{S}_{2n}^{++}$
 \Big( $\left[\begin{array}{cc}
           \bar{P}^i\!\!&\!\!\bf{0} \\ \bf{0}\!\!&\!\!P^i\end{array}\right]
           \in \mathscr{S}_{2n}^{+}$  \Big)
 under exact observability (exact detectability). \\
 Evaluating the Lyapunov recursion (\ref{pe2})  at
$i=i+1$ and subtracting equation (\ref{pe21}),   the following recursion is
derived:
\begin{eqnarray} \label{pe22}
 && (\mathcal{A}_1^{F^{i+1},\bar{F}^{i+1}})'
              \left[\begin{array}{cc}
         \!   \Delta \bar{P}^i  \!\!&\!\! \bf{0} \\ \!\bf{0}\!\!   &\!\!  \Delta P^i\!\end{array}\right]\!(\mathcal{A}_1^{F^{i+1},\bar{F}^{i+1}})
          +(\mathcal{A}_2^{F^{i+1}})'\nonumber\\
            &&
            \cdot \left[\begin{array}{cc}
          \!   \Delta \bar{P}^i\!\! &\!\!   \bf{0}\! \\\! \bf{0} \!\!  &\!\!\Delta P^i\!\end{array}\right]\! (\mathcal{A}_2^{F^{i+1}})
            \!+\!(\mathcal{A}_3^{F^{i+1},\bar{F}^{i+1}})' \! \left[\begin{array}{cc}\!
            \Delta \bar{P}^i \!\! &\!\!  \bf{0}\! \\ \!\bf{0}  \!\!  & \!\!\Delta  P^i\!\end{array}\right]\nonumber\\
            &&
            \cdot
             (\mathcal{A}_3^{F^{i+1},\bar{F}^{i+1}})\!- \!\left[\begin{array}{cc}
            \!\Delta \bar{P}^i\!\!  &\!\!  \bf{0}\! \\\! \bf{0} \!\! & \!\! \Delta P^i\!\end{array}\right]\nonumber\\
            &&
            =-\left[\begin{array}{cc}
          \!   \Delta \hat{F}^{i} \!\! &\!\!  \bf{0}\! \\\! \bf{0}  \!\!& \!\!\Delta F^{i}\!\end{array}\right]'\!
              \left[\begin{array}{cc}
           \! \Upsilon^{P^i, \bar{P}^i} \!\!  & \!\! \bf{0} \!\\ \!\bf{0}\!\!   & \!\! \Upsilon^{P^i}\!\end{array}\right]\!
            \left[\begin{array}{cc}
             \!\Delta \hat{F}^{i} \!\! & \!\!\bf{0}\!\\ \!\bf{0}\!\!  &\!\! \Delta F^{i}\!\end{array}\right],
\end{eqnarray}
where
$ \Delta \bar{P}^i =\bar{P}^{i}\!-\!\bar{P}^{i+1}$,
$ \Delta  {P}^i = {P}^{i}\!- \!{P}^{i+1}$,
$ \Delta \bar{F}^i =\bar{F}^{i}\!-\!\bar{F}^{i+1}$, and
$ \Delta  {F}^i = {F}^{i}\!-\! {F}^{i+1}$.
  {Based on Theorem 4 in \cite{chenchen}}, there exists a unique solution
 $\left[\begin{array}{cc}
        \!   \Delta  \bar{P}^{i}  \!\!&\!\!  \bf{0} \!\\ \!\bf{0}  \!\!&  \!\! \Delta  P^{i}\!\end{array}\right]\!\in\mathscr{S}_{2n}^{++}$
            \Big($\left[\begin{array}{cc}
          \!  \Delta \bar{P}^{i} \!\! & \!\!\bf{0} \!\\\! \bf{0} \!\! &\!\!   \Delta  P^{i}\!\end{array}\right]\!\in\mathscr{S}_{2n}^{+}$\Big) to
             (\ref{pe22}), which results in that $\{(P^i, \bar{P}^i)\}_{i=0}^
             \infty$ is a monotonically decreasing sequence.
             Moreover, the sequence $\{(P^i, \bar{P}^i)\}_{i=0}^
             \infty$ is with low bound zero, it must have a unique limit
             $(P^*, \bar{P}^*)$ satisfying GAREs (\ref{gare1})-(\ref{gare2}),
             i.e., $\lim_{i\rightarrow \infty}P^{i}=P^*\leq P^{i+1}\leq P^i$
             and $\lim_{i\rightarrow \infty}\bar{P}^{i}=\bar{P}^*
             \leq \bar{P}^{i+1}\leq \bar{P}^i$. \\
             Taking the limit on both sides of (\ref{pim1})-(\ref{pim2}), we have
             \begin{eqnarray*}
             \lim_{i\rightarrow \infty} F^{i+1}
             &=&\lim_{i\rightarrow \infty}-[\Upsilon^{P^i}]^{-1}[M^{P^i}]\\
             &
             =&\lim_{i\rightarrow \infty}-[\Upsilon^{P^*}]^{-1}[M^{P^*}]
              :=F^*
             \end{eqnarray*}
            and
            \begin{eqnarray*}
             \lim_{i\rightarrow \infty} \bar{F}^{i+1}
             &\!=&\!\lim_{i\rightarrow \infty}\!-\!\{[\Upsilon^{P^i, \bar{P}^i}]^{-1}\!
             [M^{P^i, \bar{P}^i}]\!
-\![\Upsilon^{P^i}]^{-1}\![M^{P^i}]\}\\
             & 
             \!\!=\!\!& \lim_{i\!\rightarrow \!\infty}\!-\!\{[\Upsilon^{P^*, \bar{P}^*}]^{\!-1}
             [M^{P^*, \bar{P}^*}]
\!-\![\Upsilon^{P^*}]^{\!-1}\![M^{P^*}]\}\\
&\!\! :=\!\!&\bar F^*,
             \end{eqnarray*}
 where $F^*$ and $\bar{F}^*$ satisfy (\ref{gain1}) and (\ref{gain2}), respectively.
\hfill $\square$
\section*{ {Appendix B}}
 {\textit{Proof of   Proposition \ref{proposition123}.}}\\
  Based on the relationship between the  solution of  generalized Lyapunov equation and the ASMS discussed in Theorem 4 of \cite{chenchen}, it is sufficient to prove that the following pair is exactly detectable:
 $$[(\mathscr{A}_1^{F, \bar{F}})';(\mathscr{A}_2^{F})', (\mathscr{A}_3^{F})'\mid (\mathbb{F}_1
Z_2\mathbb{F}_1'
 +\mathbb{F}_2(Z_1-Z_2)\mathbb{F}_2')].
 $$
For Primal Problem  I, we  seek a  matrix $\tilde{\mathbbm{S}} $ satisfying:
 \begin{eqnarray}
 &
 (\mathscr{A}_1^{F, \bar{F}})'\!\tilde{\mathbbm{S}}\mathscr{A}_1^{F, \bar{F}}\!+\! (\mathscr{A}_2^{F})'\tilde{\mathbbm{S}}
 \mathscr{A}_2^{F} \! +\! (\mathscr{A}_3^{F})'\tilde{\mathbbm{S}} \mathscr{A}_3^{F}\!=\!\lambda \!\tilde{\mathbbm{S}},\label{phb11}\\
 &
  [\mathbb{F}_1
Z_2\mathbb{F}_1'
 +\mathbb{F}_2(Z_1-Z_2)\mathbb{F}_2']\tilde{\mathbbm{S}}=0 \label{phb22}
 \end{eqnarray}
 with $\lambda\geq |1|$.
Given  the block structure of   $\tilde{\mathbbm{S}}$ matches that of
$\tilde{\mathbbm{S}}_{\mathcal{P}_I}$,   equation (\ref{phb22}) leads to
\begin{eqnarray*}
\begin{array}{l}
\left[\begin{array}{cccc}
Z_2(\bar{\mathbbm{S}}^{11}+ \hat{F}'\bar{\mathbbm{S}}^{12'})
&
Z_2(\bar{\mathbbm{S}}^{12}+ \hat{F}'\bar{\mathbbm{S}}^{22})\\
\hat{F}Z_2(\bar{\mathbbm{S}}^{11}+
 \hat{F}'\bar{\mathbbm{S}}^{12'})&
\hat{F}Z_2(\bar{\mathbbm{S}}^{12}
+ \hat{F}'\bar{\mathbbm{S}}^{22})
\end{array}
\right]=0,\\
\left[\begin{array}{cccc}
\!(Z_1\!-\!Z_2)\!({\mathbbm{S}}^{11}\!+\!{F}' {\mathbbm{S}}^{12'})
\!&\!
(Z_1\!-\!Z_2)\!( {\mathbbm{S}}^{12}\!+\!  {F}' {\mathbbm{S}}^{22})\\
 \!{F}(Z_1\!-\!Z_2)\!({\mathbbm{S}}^{11}\!+\!
 {F}' {\mathbbm{S}}^{12'})\!&\!
 {F}(Z_1\!-\!Z_2)({\mathbbm{S}}^{12}
\!+ \!{F}' {\mathbbm{S}}^{22})
\end{array}\!\!
\right]=0.
\end{array}
 \end{eqnarray*}
Since $Z_1>Z_2>0$,   it follows that:
 \begin{eqnarray*}
 &&\bar{\mathbbm{S}}^{11}+ \hat{F}'\bar{\mathbbm{S}}^{12'}=0,\
   \bar{\mathbbm{S}}^{12}+ \hat{F}'\bar{\mathbbm{S}}^{22}=0, \\
 && {\mathbbm{S}}^{11}+{F}' {\mathbbm{S}}^{12'}=0, \
 {\mathbbm{S}}^{12}+  {F}' {\mathbbm{S}}^{22}=0.
 \end{eqnarray*}
Through straightforward calculations, we obtain:
 $$(\mathscr{A}_1^{F, \bar{F}})'\tilde{\mathbbm{S}}\mathscr{A}_1^{F, \bar{F}}=(\mathscr{A}_2^{F})'\tilde{\mathbbm{S}}
 \mathscr{A}_2^{F} = (\mathscr{A}_3^{F})'\tilde{\mathbbm{S}} \mathscr{A}_3^{F}=0,$$
 which implies that  the solution is  $\mathbbm{S}=0$.
 Using  the stochastic PBH eigenvector test, we confirm the required exact detectability in Primal Problem  I. According to  {the assertions (1) and (4) in Theorem 4 of \cite{chenchen}}, $(F, \bar{F})\in\mathcal{F}$  iff the GLE (\ref{gle1})  has  a unique solution $\mathbbm{S}\geq 0$. \hfill $\square$

\section*{ {Appendix C}}
 {\textit{Proof of   Proposition \ref{proposition1}.}}\\
 The proof  consists of two parts: one for Primal Problem  I and the other for Primal Problem  II.\\
Part A: For Primal Problem  I.
  Set
\begin{eqnarray}\label{sdefinition}
\tilde{\mathbbm{S}}\triangleq \sum_{l=1}^r\sum_{k=\psi}^\infty
\mathbb{E}[(\mathbb{V}_k^{F, \bar{F}; \mathbb{V}_\psi^{z_l}})
(\mathbb{V}_k^{F, \bar{F}; \mathbb{V}_\psi^{z_l}})'].
\end{eqnarray}
Clearly, $\tilde{\mathbbm{S}}$ is  a block digaram matrix with the form
\begin{eqnarray*}
 \tilde{\mathbbm{S}}&\triangleq&
\left[\begin{array}{ccccc}
\bar{\mathbbm{S}}& {\bf 0}\\
{\bf 0}& {\mathbbm{S}} \end{array}\right]_{\bar{\mathbbm{S}}\in\mathscr{S}_{n+m}, {\mathbbm{S}}\in\mathscr{S}_{n+m}}\\
 &
=&\sum_{l=1}^r\sum_{k=\psi}^\infty
\mathbb{E}
\left[\begin{array}{ccccc}
\!\mathbb{E}v_k(\mathbb{E}v_k)'\!\!&\!\! {\bf 0}\!\\
\!{\bf 0}\!\!&\!\!(v_k-\mathbb{E}v_k)\!(v_k-\mathbb{E}v_k)' \! \end{array}\right]
\end{eqnarray*}
Along with the system (\ref{ausys}), there is
\begin{eqnarray*}
 \tilde{\mathbbm{S}}&=&
\mathbb{F}_1
\sum_{l=1}^r [\mathbb{E}z_l(\mathbb{E}z_l)']\mathbb{F}_1' + \sum_{l=1}^r\sum_{k=\psi}^\infty
\mathbb{E}[(\mathbb{V}_{k+1}^{F, \bar{F}; \mathbb{V}_\psi^{z_l}})\\
& 
\cdot & (\mathbb{V}_{k+1}^{F, \bar{F}; \mathbb{V}_\psi^{z_l}})']
+\mathbb{F}_2
\sum_{l=1}^r \mathbb{E}[(z_l-\mathbb{E}z_l)(z_l-\mathbb{E}z_l)']
\mathbb{F}_2'.
\end{eqnarray*}
According to
\begin{eqnarray*}
&&\mathbb{E}[(\mathbb{V}_{k+1}^{F, \bar{F}; \mathbb{V}_\psi^{z_l}})
(\mathbb{V}_{k+1}^{F, \bar{F}; \mathbb{V}_\psi^{z_l}})']\\
&&
= \mathscr{A}_1^{F, \bar{F}}\mathbb{E}[(\mathbb{V}_{k}^{F, \bar{F}; \mathbb{V}_\psi^{z_l}})
(\mathbb{V}_{k}^{F, \bar{F}; \mathbb{V}_\psi^{z_l}})']
(\mathscr{A}_1^{F, \bar{F}})'\\
&& \  \ \ \ \ \
+ \mathscr{A}_2^{F}\mathbb{E}[(\mathbb{V}_{k}^{F, \bar{F}; \mathbb{V}_\psi^{z_l}})
(\mathbb{V}_{k}^{F, \bar{F}; \mathbb{V}_\psi^{z_l}})']
(\mathscr{A}_2^{F})'\\
&&\  \ \ \ \ \
+ \mathscr{A}_3^{F}\mathbb{E}[(\mathbb{V}_{k}^{F, \bar{F}; \mathbb{V}_\psi^{z_l}})
(\mathbb{V}_{k}^{F, \bar{F}; \mathbb{V}_\psi^{z_l}})']
(\mathscr{A}_3^{F})',
\end{eqnarray*}
the constraint condition GLE (\ref{gle1}) is obtained.
Due to $(F, \bar{F})\in \mathcal{F}$, 
there is  a unique solution $0\leq \tilde{\mathbbm{S}}<\infty$, which
coincides with the definition of matrix
$\tilde{\mathbbm{S}}$  given in (\ref{sdefinition}).  The objective function in Problem 1 can be equivalently expressed as
\begin{eqnarray*}
&&\hat{J}(z_1, \cdots, z_r; F, \bar{F})\\
&& =\mathbb{E}
\sum_{l=1}^r
\sum\limits_{k=\psi}^\infty
(\mathbb{V}_k^{F, \bar{F}; \mathbb{V}_\psi^{z_l}})'
\left[\begin{array}{ccccc}
\!\Lambda\!+\!\bar{\Lambda}\!\!&\!\!\bf{0}\!\\\!
\bf{0}\!\! &\!\!\Lambda\! \end{array}\right]
\mathbb{V}_k^{F, \bar{F}; \mathbb{V}_\psi^{z_l}}\\
&&
=Tr\Big(\left[\begin{array}{ccccc}
\!\Lambda\!+\!\bar{\Lambda}\!\!&\!\!\bf{0} \!\\\!
\bf{0} \!\!&\!\!\Lambda \!\end{array}\right]\!\tilde{\mathbbm{S}}\!\Big).
\end{eqnarray*}
From the above discussion,  one can directly obtained that  $J_{\mathcal{P}_I}=\hat{J}(F^*, \bar{F}^*)=Tr\Big(\left[\begin{array}{ccccc}
\!\Lambda\!+\!\bar{\Lambda}\!\!& \!\!\bf{0} \!\\
 \!\bf{0} \!\!&\!\!\Lambda\! \end{array}\right]\!\tilde{\mathbbm{S}}_{\mathcal{P}_I}^*\Big)$, which
 ($\tilde{\mathbbm{S}}_{\mathcal{P}_I}^*, F^*, \bar{F}^*$)
is the unique solution of  the GLE (\ref{gle1}), with $(F^*, \bar{F}^*)$ being
given as (\ref{gain1})-(\ref{gain2}).\\
Part B: For Primal Problem  II.
For condition $(F, \bar{F})\in\mathcal{F}$ holding, it is directly know that
$[\mathscr{A}_1^{F, \bar{F}};\mathscr{A}_2^{F}, \mathscr{A}_3^{F}]$ is ASMS. Combing with $\left[\begin{array}{ccccc}
\!\Lambda\!+\!\bar{\Lambda}\!\!& \!\!\bf{0} \!\\\!
0\!\!&\!\!\Lambda \!\end{array}\right]\geq  \bf{0} $,
 the GLE (\ref{gle2}) admits a unique solution
$\tilde{\mathscr{X}}\in\mathscr{S}_{2n+2m}^+$.  By the completing square method,
we have
\begin{eqnarray*}
&&\hat{J}(z_1, \cdots, z_r; F, \bar{F})\\
&& =\!\mathbb{E}\!
\sum_{l=1}^r
\sum\limits_{k=\psi}^\infty
\Big[(\mathbb{V}_k^{F, \bar{F}; \mathbb{V}_\psi^{z_l}})'
\left[\begin{array}{ccccc}
\!\Lambda\!+\!\bar{\Lambda}\!\!& \!\!\bf{0} \!\\
\!\bf{0} \!\!&\!\!\Lambda\! \end{array}\!\right]\!
\mathbb{V}_k^{F, \bar{F}; \mathbb{V}_\psi^{z_l}}\!+\!\Delta V_k^{F, \bar{F}, z_l}\!\Big]\\
&& +\sum_{l=1}^r\Big\{\mathbb{E}\big[ \!(\mathbb{V}_\psi^{z_l})'\!\tilde{\mathscr{X}} \! \mathbb{V}_\psi^{z_l}\big]\!-\!\lim\limits_{k\rightarrow \infty}\mathbb{E}\![ (\mathbb{V}_k^{F, \bar{F}; \mathbb{V}_\psi^{z_l}})'\!\tilde{\mathscr{X}} \!\mathbb{V}_k^{F, \bar{F}; \mathbb{V}_\psi^{z_l}}\!\big]\!\Big\}\\
&&=
\sum_{l=1}^r
\sum\limits_{k=\psi}^\infty \mathbb{E}\big[(\mathbb{V}_k^{F, \bar{F}; \mathbb{V}_\psi^{z_l}})' \!\Pi\! \mathbb{V}_k^{F, \bar{F}; \mathbb{V}_\psi^{z_l}}\!\big]\!+\!\sum_{l=1}^r \mathbb{E}\!\big[ \! (\mathbb{V}_\psi^{z_l})'\!\tilde{\!\mathscr{X}}\! \mathbb{V}_\psi^{z_l}\!]\\
&&=\sum_{l=1}^r \mathbb{E}\big[ (\mathbb{V}_\psi^{z_l})'\tilde{\mathscr{X}} \mathbb{V}_\psi^{z_l}],
\end{eqnarray*}
where $\Pi\!\triangleq \! (\mathscr{A}_1^{F, \bar{F}})'\!\tilde{\mathscr{X}}\mathscr{A}_1^{F, \bar{F}}
\!+\!(\mathscr{A}_2^{F})'\tilde{\mathscr{X}}\mathscr{A}_2^{F}
\!+\!(\mathscr{A}_3^{F})'\tilde{\mathscr{X}}\mathscr{A}_3^{F}\!+\!
\left[\begin{array}{ccccc}
\!\Lambda\!+\!\bar{\Lambda}\!&\! \bf{0} \!\\
 \!\bf{0} \!&\!\Lambda \!\end{array}\right]\!-\!\tilde{\mathscr{X}}\!=\!0$, $\Delta V_k^{F, \bar{F}, z_l}=\mathbb{E}\big[ (\mathbb{V}_{k+1}^{z_l})'\tilde{\mathscr{X}} \mathbb{V}_{k+1}^{z_l}\big]-\mathbb{E}\big[ (\mathbb{V}_{k}^{z_l})'
\tilde{\mathscr{X}} \mathbb{V}_{k}^{z_l}\big]$.
Since $\tilde{\mathscr{X}}$, $F$, and $\bar{F}$ satisfy GLE (\ref{gle2}), there is
\begin{eqnarray*}
&&\hat{J}(F^*, \bar{F}^*)=
\inf\limits_{X\in\mathscr{S}_{2n+2m}^+, F, \bar{F}\in\mathcal{F}}
\sum_{l=1}^r \mathbb{E}\big[ (\mathbb{V}_\psi^{z_l})'\tilde{\mathscr{X}} \mathbb{V}_\psi^{z_l}]\\
&& \inf\limits_{X\in\mathscr{S}_{2n+2m}^+, F, \bar{F}\in\mathcal{F}}
Tr (\hbar)=J_{\mathcal{P}_{II}}.
\end{eqnarray*}
\hfill $\square$

 \section*{ {Appendix D}}
 {\textit{Proof of   Lemma \ref{lemmaproblem2}.}}\\
Considering  the matrix blocks $\bar{\mathbbm{S}}_{\mathcal{P}_I}$  and ${\mathbbm{S}}_{\mathcal{P}_I}$,
  the GLE  (\ref{gle1}) in Primal Problem  I can be rewritten as
\begin{eqnarray}
&&  \tilde{\mathbbm{S}}_{\mathcal{P}_I}=
 \left[\begin{array}{ccccc}
\!I_n\!\\ \!F\!+\!\bar{F} \!\end{array}\right]
 \left[\begin{array}{ccccc}
\!\hat{A}_1\!&\! \hat{B}_1 \!\end{array}\right]
\!\bar{\mathbbm{S}}_{\mathcal{P}_I}
\! \left[\begin{array}{ccccc}
\!\hat{A}_1'\!\\\! \hat{B}_1' \!\end{array}\right]
 \!\left[\begin{array}{ccccc}
\!I_n\!\\ \!F\!+\!\bar{F}\!\end{array}\right]'\nonumber\\
&& \ \ \ \  \ \ \ \
\!+\!\left[\begin{array}{ccccc}
\!I_n\!\\\! F\!+\!\bar{F}\!\end{array}\right]\!Z_2
\!\left[\begin{array}{ccccc}
\!I_n\!\\ \!F\!+\!\bar{F}\!\end{array}\right]'\label{gle12}
  \end{eqnarray}
and
\begin{eqnarray}
&& {\mathbbm{S}}_{\mathcal{P}_I}=
 \left[\begin{array}{ccccc}
\!I_n\!\\\! F \!\end{array}\right]
\! \left[\begin{array}{ccccc}
\!{A}_1\!&\! {B}_1 \!\end{array}\right]
\!{\mathbbm{S}}_{\mathcal{P}_I}
\! \left[\begin{array}{ccccc}
\!{A}_1'\!\\ \!{B}_1' \!\end{array}\right]
 \!\left[\begin{array}{ccccc}
\!I_n\!\\\! F\!\end{array}\right]'\nonumber\\
&& \ \ \ \  \ \ \ \
+ \!\left[\begin{array}{ccccc}
\!I_n\!\\\! F\! \end{array}\right]
 \!\left[\begin{array}{ccccc}
\!{A}_2\!&\! {B}_2 \!\end{array}\right]
\!{\mathbbm{S}}_{\mathcal{P}_I}
 \!\left[\begin{array}{ccccc}
\!{A}_2'\!\\\! {B}_2' \!\end{array}\right]
\! \left[\begin{array}{ccccc}
\!I_n\!\\\! F\!\end{array}\right]'\nonumber\\
&& \ \ \ \  \ \ \ \
\!+\! \left[\begin{array}{ccccc}
\!I_n\!\\\! F\! \end{array}\right]
\! \left[\begin{array}{ccccc}
\!\hat{A}_2\!&\! \hat{B}_2 \!\end{array}\right]
\!\bar{\mathbbm{S}}_{\mathcal{P}_I}
 \!\left[\begin{array}{ccccc}
\!\hat{A}_2'\!\\ \!\hat{B}_2' \!\end{array}\right]
\! \left[\begin{array}{ccccc}
\!I_n\!\\ \!F\!\end{array}\right]'\nonumber\\
&& \ \ \ \  \ \ \ \
\!+\!\left[\begin{array}{ccccc}
\!I_n\!\\\! F \!\end{array}\right](Z_1\!-\!Z_2)
\left[\begin{array}{ccccc}
\!I_n\!\\ \!F\! \end{array}\right]'.\label{gle13}
  \end{eqnarray}
By comparing the first $n\times n$ sub-block matrix of the above two equatons,
there are
\begin{eqnarray}\label{barss}
\left\{\begin{array}{l}
\bar{\mathbbm{S}}_{\mathcal{P}_I}^{11}=
\left[\begin{array}{ccccc}
\!\hat{A}_1\!\!& \!\!\hat{B}_1 \!\end{array}\right]
\!\bar{\mathbbm{S}}_{\mathcal{P}_I}
 \!\left[\begin{array}{ccccc}
\!\hat{A}_1'\!\!\\\!\! \hat{B}_1' \!\end{array}\right]\!+\!Z_2,\\
{\mathbbm{S}}_{\mathcal{P}_I}^{11}=
\! \left[\begin{array}{ccccc}
\!{A}_1\!\!&\!\! {B}_1 \!\end{array}\right]
\!{\mathbbm{S}}_{\mathcal{P}_I}
\! \left[\begin{array}{ccccc}
\!{A}_1'\!\!\\ \!\!{B}_1' \!\end{array}\right]
+
\! \left[\begin{array}{ccccc}
\!{A}_2\!\!&\!\! {B}_2\! \end{array}\right]
\!{\mathbbm{S}}_{\mathcal{P}_I}
 \!\left[\begin{array}{ccccc}
\!{A}_2'\!\!\\ \!\!{B}_2' \!\end{array}\right]
\\
\ \ \ \ \ \ \ \ \ \ \ +
 \!\left[\begin{array}{ccccc}
\!\hat{A}_2\!\!& \!\!\hat{B}_2\! \end{array}\right]
\!\bar{\mathbbm{S}}_{\mathcal{P}_I}
\! \left[\begin{array}{ccccc}
\!\hat{A}_2'\!\!\\ \!\!\hat{B}_2' \!\end{array}\right]
\!+\!
(Z_1\!-\!Z_2).
\end{array}
\right.
\end{eqnarray}
Using the rewritten version of GLE (\ref{gle1}) given in (\ref{gle12})-(\ref{gle13}),
we have
\begin{eqnarray*}
&&\bar{\mathbbm{S}}_{\mathcal{P}_I}^{11}\\
&&=
\left[\begin{array}{ccccc}
\!\hat{A}_1\!\!&\!\! \hat{B}_1\! \end{array}\right]
\!\Big\{\!
\left[\begin{array}{ccccc}
\!I_n\!\\ \!F\!+\!\bar{F}\! \end{array}\right]
 \!\left[\begin{array}{ccccc}
\!\hat{A}_1\!\!& \!\!\hat{B}_1\! \end{array}\right]
\!\bar{\mathbbm{S}}_{\mathcal{P}_I}
 \!\left[\begin{array}{ccccc}
\!\hat{A}_1'\!\\\! \hat{B}_1' \!\end{array}\right]\\
&&\ \ \ \cdot
 \!\left[\begin{array}{ccccc}
\!I_n\!\\ \!F\!+\!\bar{F}\!\end{array}\right]'
\!+\!\left[\begin{array}{ccccc}
\!I_n\!\\ \!F\!+\!\bar{F}\!\end{array}\right]\!Z_2\!
\!\left[\begin{array}{ccccc}
\!I_n\!\\\! F\!+\!\bar{F}\!\end{array}\right]'
\!\Big\}\!\left[\begin{array}{ccccc}
\!\hat{A}_1'\!\\\! \hat{B}_1'\! \end{array}\right]\!+\!Z_2\\
&&
=[\hat{A}_1\!+\!\hat{B}_1\!(F\!+\!\bar{F})]\!\bar{\mathbbm{S}}_{\mathcal{P}_I}^{11}
\![\hat{A}_1\!+\!\hat{B}_1(F+\bar{F})]'\!+\!Z_2.
\end{eqnarray*}
Similarly,
\begin{eqnarray*}
&&{\mathbbm{S}}_{\mathcal{P}_I}^{11}
=[{A}_1+B_1F]{\mathbbm{S}}_{\mathcal{P}_I}^{11}
[{A}_1+B_1F]'+[{A}_2+B_2F]\\
&&\ \ \ \ \  \ \ \ \ \ \cdot {\mathbbm{S}}_{\mathcal{P}_I}^{11}
[{A}_2+B_2F]'
+[\hat{A}_2+\hat{B}_2(F+\bar{F})]\bar{\mathbbm{S}}_{\mathcal{P}_I}^{11}\\
&&\ \ \ \ \  \ \ \ \ \ \cdot
[\hat{A}_2+\hat{B}_2(F+\bar{F})]'
+Z_1-Z_2.
\end{eqnarray*}
Combing the notations $\mathbb{A}_1$, $\mathbb{A}_2$,
and $\mathbb{A}_3$, the
first result is proved. Next, for the equation representations of
$\bar{\mathbbm{S}}_{\mathcal{P}_I}^{11}$ and ${\mathbbm{S}}_{\mathcal{P}_I}^{11}$ in (\ref{barss}),
there are
\begin{eqnarray*}
\bar{\mathbbm{S}}_{\mathcal{P}_I}&=&
\left[\begin{array}{ccccc}
I_n\\ F+\bar{F}\end{array}\right]\bar{\mathbbm{S}}_{\mathcal{P}_I}^{11}
\left[\begin{array}{ccccc}
I_n\\ F+\bar{F}\end{array}\right]'\\
&=&
\left[\begin{array}{ccccc}
\bar{\mathbbm{S}}_{\mathcal{P}_I}^{11}&
\bar{\mathbbm{S}}_{\mathcal{P}_I}^{11}( F+\bar{F})'\\ *&
( F+\bar{F})\bar{\mathbbm{S}}_{\mathcal{P}_I}^{11}( F+\bar{F})'\end{array}\right],\\
{\mathbbm{S}}_{\mathcal{P}_I}&=&
\left[\begin{array}{ccccc}
I_n\\ F\end{array}\right]{\mathbbm{S}}_{\mathcal{P}_I}^{11}
\left[\begin{array}{ccccc}
I_n\\ F\end{array}\right]'\\
&=&
\left[\begin{array}{ccccc}
{\mathbbm{S}}_{\mathcal{P}_I}^{11}&
{\mathbbm{S}}_{\mathcal{P}_I}^{11}F'\\ *&
 F{\mathbbm{S}}_{\mathcal{P}_I}^{11}F'\end{array}\right],
\end{eqnarray*}
which result  in $\bar{\mathbbm{S}}_{\mathcal{P}_I}^{11}( F+\bar{F})'=
\bar{\mathbbm{S}}_{\mathcal{P}_I}^{12}$ and ${\mathbbm{S}}_{\mathcal{P}_I}^{11}F'={\mathbbm{S}}_{\mathcal{P}_I}^{12}$.
Since $Z_1>Z_2>0$,  the results (\ref{fsfs1}) and  (\ref{fsfs2}) are obtained.
\hfill $\square$
\section*{ {Appendix E}}
 {\textit{Proof of   Theorem   \ref{sdth}.}}\\
 Strong duality is established by demonstrating that the duality gap for Primal Problem
 I is zero.
Note that the dual  function $d(\tilde{\mathfrak{X}})$ can be rewritten as
\begin{eqnarray*}
&&d(\tilde{\mathfrak{X}})=\inf\limits_{\tilde{\mathbbm{S}}\in\mathscr{S}_{2n+2m}^+, F, \bar{F}\in\mathcal{F}}Tr [Z_2\mathbb{F}_1'\tilde{\mathfrak{X}}\mathbb{F}_1
+
(Z_1-Z_2)\mathbb{F}_2'\tilde{\mathfrak{X}}\mathbb{F}_2]\\
&& \ \ \ \  \ \
+Tr\big\{[ (\mathscr{A}_1^{F, \bar{F}})'\tilde{\mathfrak{X}}  \mathscr{A}_1^{F, \bar{F}}+(\mathscr{A}_2^{F})'\tilde{\mathfrak{X}}
\mathscr{A}_2^{F}
 +(\mathscr{A}_3^{F})'\tilde{\mathfrak{X}}
\mathscr{A}_3^{F}\\
&& \ \ \ \ \ \  +\left[\begin{array}{ccccc}
\Lambda+\bar{\Lambda}&\bf{0}\\
\bf{0} &\Lambda \end{array}\right]-\tilde{\mathfrak{X}}]\tilde{\mathbbm{S}}\big\}
\\
&& \ \ \ \ \ \
=\left\{\begin{array}{l}
\inf\limits_{F, \bar{F}\in\mathcal{F}}Tr [Z_2\mathbb{F}_1'\tilde{\mathfrak{X}}\mathbb{F}_1
+
(Z_1-Z_2)\mathbb{F}_2'\tilde{\mathfrak{X}}\mathbb{F}_2], \tilde{\mathfrak{X}}\in \digamma,
\\
-\infty, \ \ \ \ \ otherwise
\end{array}
\right.
\end{eqnarray*}
with
\begin{eqnarray*}
&& \digamma\triangleq \big\{
\tilde{\mathfrak{X}}\in
\mathscr{S}_{2n+2m}^+:   (\mathscr{A}_1^{F, \bar{F}})'\tilde{\mathfrak{X}}  \mathscr{A}_1^{F, \bar{F}}\!+\!(\mathscr{A}_2^{F})'\tilde{\mathfrak{X}}
\mathscr{A}_2^{F}\\
&& \ \ \ \ \  
 \!+\!(\mathscr{A}_3^{F})'\tilde{\mathfrak{X}}
\mathscr{A}_3^{F} \! +\!\left[\begin{array}{ccccc}
\!\Lambda\!+\!\bar{\Lambda}\!\!&\!\! \bf{0}\! \\
\bf{0} \!\!&\!\!\Lambda \!\end{array}\right]\!-\!\tilde{\mathfrak{X}}\!\geq\! 0, \forall F, \bar{F}\in\mathcal{F}
\big\}.
\end{eqnarray*}
Similar analysis as in \cite{lee2019,liman}, we know that the admissible set $\digamma$ is nonempty with $\tilde{\mathfrak{X}}_{\mathcal{P}_{II}}^*\in\digamma$. Hence, the objective value of Dual Problem is changed as
\begin{eqnarray*}
J_{\mathcal{D}}&=&\sup\limits_{\tilde{\mathfrak{X}}\in\digamma}
d(\tilde{\mathfrak{X}})\\
&=&\sup\limits_{\tilde{\mathfrak{X}}\in\digamma}
\inf\limits_{F, \bar{F}\in\mathcal{F}}Tr [Z_2\mathbb{F}_1'\tilde{\mathfrak{X}}\mathbb{F}_1
+
(Z_1-Z_2)\mathbb{F}_2'\tilde{\mathfrak{X}}\mathbb{F}_2] .
\end{eqnarray*}
 Since now, $d(\tilde{\mathfrak{X}})$ is a quadratic function with respect to $F, \bar{F}\in\mathcal{F}$
Clearly,
$J_{\mathcal{P}_{II}}=d(\tilde{\mathfrak{X}}_{\mathcal{P}_{II}}^*)\leq J_{\mathcal{D}}$.  Combing with the condition $J_{\mathcal{P}_{I}}=J_{\mathcal{P}_{II}}$ proposed  in Proposition \ref{proposition1},
the strong duality $J_{\mathcal{P}_{I}}=J_{\mathcal{D}}$ holds.\hfill $\square$
\section*{ {Appendix F}}
 {\textit{Proof of   Proposition \ref{prpddddd}.}}\\
For the Problem 5, the KKT conditions can be summarized as
\begin{description}
  \item[(i)] the primal feasibility condition given as (\ref{kkt1})-(\ref{kkt2}).
  \item[(ii)]  the complementary slackness condition $Tr(\tilde{\mathbbm{S}}\tilde{\mathfrak{X}}_0)=0$, in where the Lagrange multiplier $\tilde{\mathfrak{X}}_0=0$ since $\tilde{\mathbbm{S}}>0$.
  \item[(iii)] the Lagrangian gradient condition
  \begin{align*}
 &\frac{\partial \hat{L}(\tilde{\mathfrak{X}}, \tilde{\mathfrak{X}}_0, F, \bar{F}, \tilde{\mathbbm{S}})}{\partial \tilde{\mathbbm{S}}} \Big| _{\tilde{\mathbbm{S}}=\tilde{\mathbbm{S}}^*}\\
  &
  =(\mathscr{A}_1^{F, \bar{F}})'\tilde{\mathfrak{X}}
 \mathscr{A}_1^{F, \bar{F}}
+(\mathscr{A}_2^{F})'\tilde{\mathfrak{X}}
\mathscr{A}_2^{F}
+(\mathscr{A}_3^{F})'\tilde{\mathfrak{X}}
\mathscr{A}_3^{F}\\
& \ \ \ \ +\left[\begin{array}{ccccc}
\!\Lambda\!+\!\bar{\Lambda}\!\!&\!\! \bf{0}\!\\
\!\bf{0} \!\!&\!\!\Lambda \!\end{array}\right]\!-\!\tilde{\mathfrak{X}}
 \!-\!\tilde{\mathfrak{X}}_0\!=\!0,\\
 &
 \frac{\partial \hat{L}(\tilde{\mathfrak{X}}, \tilde{\mathfrak{X}}_0, F, \bar{F}, \tilde{\mathbbm{S}})}{\partial \bar{F}}\Big|_{ {F}= {F}^*, \bar{F}=\bar{F}^*}\\
 & = \frac{\partial Tr\!\left\{\!
 \left[\!\begin{array}{ccc}
 \!I\!\\\!F\!+\!\bar{F}\!\end{array}\!\right]\!\!
 \left[\begin{array}{ccc}\!
 \hat{A}_1\!\!&\!\!\hat{B}_1\!\end{array}\right]\!\bar{\mathbbm{S}}\!
  \left[\begin{array}{ccc}\!
 \hat{A}_1'\!\\\! \hat{B}_1'\!\end{array}\right]\!
 \left[\begin{array}{ccc}
\! I\!\!&\!\!(F+\bar{F})'\!\end{array}\right]\!\bar{\mathfrak{X}}\!\right\}}
 {\partial \bar{F}} \\
 &= 2[(\bar{\mathfrak{X}}_{12}^*)'\!+\!\bar{\mathfrak{X}}_{22}^*
(F^*\!+\!\bar{F}^*)]\!\left[\begin{array}{ccc}\!\hat{A}_1\!\!&\!\!\hat{B}_1\!
\end{array}
\right]\!\bar{\mathbbm{S}}\! \left[\begin{array}{ccc}\!\hat{A}_1\!\!&\!\!\hat{B}_1\end{array}
\right]'\!=\!0,\\
   &  \frac{\partial \hat{L}(\tilde{\mathfrak{X}}, \tilde{\mathfrak{X}}_0, F, \bar{F}, \tilde{\mathbbm{S}})}{\partial  {F}}\Big|_{ {F}= {F}^*, \bar{F}=\bar{F}^*}\\
 & =
    \frac{\!\partial \!Tr\!\left\{\!
 \left[\begin{array}{ccc}\!
 \!I\!\!\\\!\!F\!+\!\bar{F}\!\!\end{array}\right]\!\!
 \left[\begin{array}{ccc}
 \!\hat{A}_1\!\!&\!\!\hat{B}_1\!\end{array}\right]\!\bar{\mathbbm{S}}\!
  \left[\begin{array}{ccc}
 \!\hat{A}_1'\!\!\\\!\! \hat{B}_1'\!\end{array}\right]\!
 \left[\begin{array}{ccc}
 \!I\!\!&\!\!(F+\bar{F})'\!\!\end{array}\right]\!\bar{\mathfrak{X}}\!\right\}}
 {\partial  {F}}\\
 & \ \ \  +  \frac{\partial Tr\left\{
 \left[\begin{array}{ccc}
 I\\F \end{array}\right] \Psi
 \left[\begin{array}{ccc}
 I&F'\end{array}\right] {\mathfrak{X}}\right\}}
  {\partial  {F}}\\
  & =
  2[({\mathfrak{X}}_{12}^*)'+{\mathfrak{X}}_{22}^*
F^*]\Psi=0
  \end{align*}
\end{description}
Then, the KKT conditions (\ref{kkt1})-(\ref{kkt5}) are derived.\hfill $\square$

\section*{ {Appendix G}}
 {\textit{Proof of    Theorem \ref{theorem41}.}}\\
In Algorithm 2,  considering the matrix block $\tilde{\mathfrak{X}}_i$,
there is
\begin{align}
& \left[\begin{array}{ccccc}
\!\hat{A}_1'\!\\ \!\hat{B}_1' \!\end{array}\right]\!
\bar{P}^i\!
\left[\begin{array}{ccccc}
\!\hat{A}_1\!\!& \!\!\hat{B}_1\! \end{array}\right]\!
+\!\left[\begin{array}{ccccc}
\!\hat{A}_2'\!\\ \!\hat{B}_2'\! \end{array}\right]\!
{P}^i\!
\left[\begin{array}{ccccc}
\!\hat{A}_2\!\!&\!\! \hat{B}_2\! \end{array}\right]\!
+\! \Lambda\!+\!\bar{\Lambda}\!=\!\bar{\mathfrak{X}}_i,\label{abddsss1}\\
& \left[\begin{array}{ccccc}
 \!{A}_1'\!\\  \!{B}_1' \!\end{array}\right]\!
{P}^i
\!\left[\begin{array}{ccccc}
 \!{A}_1\!\!&\!\!  {B}_1 \!\end{array}\right]
\!+\!\left[\begin{array}{ccccc}
 \!{A}_2'\!\\ \! {B}_2' \!\end{array}\right]\!
{P}^i\!
\left[\begin{array}{ccccc}
\! {A}_2\!\!&\!\!  {B}_2 \!\end{array}\right]
\!+\!\Lambda \!= \!{\mathfrak{X}}_i \label{abddsss2}
\end{align}
with
$ \bar{P}^i\!\triangleq \!\left[\begin{array}{ccccc}
\! I\!\!&\!\!  (F^i+\bar{F}^i)' \end{array}\right]\!\bar{\mathfrak{X}}_i \!\left[\begin{array}{ccccc}
 \!I\!\\\! F^i\!+\!\bar{F}^i \!\end{array}\right]$ and $  {P}^i\!\triangleq \! \left[\begin{array}{ccccc}
 \!I\!\!&\!\!  (F^i)' \end{array}\right]\! {\mathfrak{X}}_i \left[\begin{array}{ccccc}
 \!I\!\\ \!F^i \!\end{array}\right]$.
  We pre-multiply and post-multiply   equation (\ref{abddsss1}) by $\left[\begin{array}{ccccc}
 I&  (F^i+\bar{F}^i)' \end{array}\right] $
  and its transpose, respectively. Similarly, we apply $\left[\begin{array}{ccccc}
 I&  (F^i)' \end{array}\right] $
  and its transpose to   equation (\ref{abddsss2}). This gives the policy evaluation (\ref{pe2}) in Algorithm 1.
  Considering that matrices
$ {\mathfrak{X}}_i $ and $ \bar{\mathfrak{X}}_i $
are  in the block form as   (\ref{xllaaa2}),  the primal update schemes (\ref{dualupdate11})-(\ref{dualupdate22}) are identical with policy improvements
(\ref{pim1})-(\ref{pim2}). In summary, the convergence property of Algorithm 2 is equivalent with Algorithm 1. According to Theorem \ref{th22211}, the proof is directly completed. \hfill $\square$
\end{document}